\documentclass[pdflatex,sn-mathphys-num]{sn-jnl}
\usepackage{graphicx}%
\usepackage{amsmath,amssymb,amsfonts}%
\usepackage{amsthm}%
\usepackage{mathrsfs}%
\usepackage[title]{appendix}%
\usepackage{textcomp}%
\usepackage{manyfoot}%
\usepackage{float}
\usepackage{tabularx, booktabs, multirow, siunitx, colortbl}
\usepackage{algcompatible}
\usepackage[linesnumbered, lined, boxed]{algorithm2e}
\RestyleAlgo{ruled}
\SetKwComment{tcc}{/*}{*/}
\SetNoFillComment
\DontPrintSemicolon
\SetAlgoVlined  
\SetCommentSty{mycommfont}

\usepackage{subcaption}
\usepackage[compatibility]{tikztensor}
\usetikzlibrary{external}
\usepackage{pgfplots}
\usepackage{wrapfig}
\usetikzlibrary{shapes.callouts, backgrounds}
\usetikzlibrary{pgfplots.groupplots, pgfplots.statistics}

\usepgfplotslibrary{fillbetween} 
\usepackage{pgfplotstable}
\pgfplotsset{compat=1.18} 
\usepackage{comment}
\usepackage{caption}
\usepackage{subcaption}
\tikzexternalenable

\newcommand{\R}{\mathbb{R}} 
\newcommand{\C}{\mathbb{C}} 

\newcommand{\mat}[1]{\ensuremath{\mathbf{\MakeUppercase{#1}}}}  
\renewcommand{\vec}[1]{\ensuremath{\mathbf{\MakeLowercase{#1}}}}  
\newcommand{\ten}[1]{\ensuremath{\mathcal{\MakeUppercase{#1}}}} 

\newcommand{\col}[1]{\text{col}\left( #1 \right)}  
\newcommand{\rank}[1]{\text{rank}\left( #1 \right)} 

\definecolor{kulblue}{RGB}{0,85,165}
\tikzset{fancy/.style={inner color=kulblue!5!white,%
		outer color=kulblue!20!white,%
		fill faces,%
		opacity=0.95}}
\colorlet{kulblue5}{kulblue!5!white}
\colorlet{kulblue20}{kulblue!20!white}
\colorlet{kulblue50}{kulblue!50!white}
\colorlet{kulblue70}{kulblue!70!white}
\colorlet{kulblue30}{kulblue!30!black}
\colorlet{kulblue60}{kulblue!60!black}
\colorlet{kulblue90}{kulblue!90!black}

\theoremstyle{thmstyleone}%
\newtheorem{Theorem}{Theorem}
\newtheorem{Corollary}{Corollary}%
\theoremstyle{thmstyletwo}%
\newtheorem{Example}{Example}%
\newtheorem{Remark}{Remark}%
\theoremstyle{thmstylethree}%
\newtheorem{Definition}{Definition}%

\usepackage{etoolbox}

\let\oldthebibliography\thebibliography
\renewcommand{\thebibliography}[1]{%
	\oldthebibliography{#1}%
    \scriptsize               
	\setlength{\itemsep}{3.5pt}%
	\setlength{\parskip}{2pt}%
	\setlength{\parsep}{2pt}%
}

\usepackage{geometry}
\begin{document}

\title{Tensor Train Completion from Fiberwise Observations Along a Single Mode}

\author*[1,2]{\fnm{Shakir Showkat} \sur{Sofi}}\nomail \email{shakirshowkat.sofi@kuleuven.be}

\author*[1,2]{\fnm{Lieven} \sur{De Lathauwer}}\nomail \email{lieven.delathauwer@kuleuven.be}

\affil[1]{\orgdiv{Group Science, Engineering and Technology}, \orgname{KU
		Leuven Kulak}, \orgaddress{\city{Kortrijk}, \country{Belgium}}}

\affil[2]{\orgdiv{Dept. Electrical
		Engineering (ESAT)}, \orgname{KU Leuven}, \orgaddress{\city{Leuven}, \country{Belgium}}}
\unskip
\abstract{ 
 Tensor completion is an extension of matrix completion aimed at recovering a multiway data tensor by leveraging a given subset of its entries (observations) and the pattern of observation. The low-rank assumption is key in establishing a relationship between the observed and unobserved entries of the tensor. The low-rank tensor completion problem is typically solved using numerical optimization techniques, where the rank information is used either implicitly (in the rank minimization approach) or explicitly (in the error minimization approach). Current theories concerning these techniques often study probabilistic recovery guarantees under conditions such as random uniform observations and incoherence requirements. However, if an observation pattern exhibits some low-rank structure that can be exploited, more efficient  algorithms with deterministic recovery guarantees can be designed by leveraging this structure.  This work shows how to use only standard linear algebra operations to compute the tensor train decomposition of a specific type of  ``fiber-wise'' observed tensor, where some of the fibers of a tensor (along a single specific mode) are either fully observed or entirely missing, unlike the usual entry-wise observations.  From an application viewpoint, this setting is relevant when it is easier to sample or collect a multiway data tensor along a specific mode (e.g., temporal).  The proposed completion method is fast and is guaranteed to work under reasonable deterministic conditions on the observation pattern.  Through numerical experiments, we showcase  interesting applications and use cases that illustrate the effectiveness of the proposed approach\footnote[2]{This paper is an extended version of our paper published in EUSIPCO 2024~\cite{shakir2024ttfw}. See the contribution statement for the new additions.}. 
}


\maketitle

\section{Introduction} \label{sec:intro}
We live in a data-driven world, where massive amounts of data are generated and collected daily. In~addition to the huge volume and high velocity, the~structural complexity of the data is becoming so high that it renders standard techniques inadequate. Multidimensional arrays (tensors) are data structures that offer a better way to organize and analyze the data with a higher-order structure.  As~the number of dimensions of the data increases, the~memory required to store it (and the computational effort required for its analysis) increases exponentially---an obstacle known as the ``curse of dimensionality.'' Tensor decompositions offer efficient approaches to analyzing higher-order datasets, allowing for the retention of intrinsic information within the data while taming the curse of dimensionality~\cite{oseledets2009breaking, hackbusch2012tensor, grasedyck2013literature,  sidiropoulos2017tensor,   khoromskij2018tensor,  khoromskaia2018tensor, lim2021tensors, franc2022tensor,  kolda2025tensorbook, tokcan2025tensor, sofi2025bttqst}. \par 

Real-world datasets are often noisy and incomplete for various reasons, including sensor malfunctions, recording errors, constraints imposed by privacy regulations, delays in obtaining access permissions, and~deliberate incomplete sampling to meet memory and time requirements. Analyzing incomplete datasets is a big challenge, as~the missing information may affect the accuracy and reliability of the findings and, therefore, limit subsequent applications~\cite{little2019statistical}. {\em Completion} of a partially observed dataset may be thought of as a particular specification of its unobserved entries. Completion is more crucial for higher-order datasets as they are larger, increasing the chance of missing or unreliable entries. Many problems can be framed as instances of tensor completion, e.g.,~image and video inpainting, gene expression imputation, and weather and traffic data imputation (see, e.g.,~\cite{liu2013tc, song2019tc, yuan2019ttc, dropout, pan2021sclrtc, chen2020traffic}). Without~any constraints on the completed tensor, there are infinitely many ways to specify the missing entries. Therefore, to~make the estimation meaningful, it is necessary to assume that the completed tensor satisfies specific properties (e.g., low rank, minimum volume) that lower the degrees of freedom and enable a unique solution. These properties constrain the unobserved entries and help establish their relationships with the observations. The~principle of parsimony (Occam's razor) serves as a (heuristic) guideline for model selection, stating that if multiple competing models explain the same data, the~model with lower complexity is the best. The~rank can be considered a kind of natural measure of the complexity; therefore, {\em low rankness} is a reasonable perspective for promoting parsimonious models with only a few parameters explaining the~data.  \par  

The literature outlines two approaches to employing low-rank constraints in tensor completion: the rank minimization method (implicit method) and the error minimization method (explicit method). In~the former approach, rank is used implicitly as an optimization objective to minimize, with~observed entries as constraints to be satisfied. As~the rank function is non-convex and NP-hard~\cite{boyd1997sdp, fazel2004rank}, convex surrogates of the rank function are used to relax this problem, allowing for an approximate but efficient solution. In~the latter approach, a~hypothesis (tensor decomposition) model is explicitly imposed on a partially observed tensor with a specific, fixed low rank. The~loss function---typically continuous and differentiable---is minimized to fit the model parameters (see, e.g.,~\cite{candes2010mcvp, candes2012mc, liu2013tc, kressner2013lrtc, song2019tc, yuan2019ttc, gandy2011tc}). In~this case, the~flexibility in choosing a differentiable objective function enables the use of gradient-based optimization approaches (e.g., first- and/or second-order methods). Moreover, this approach explicitly uses rank information, which can be meaningful in applications where the rank has a physical significance; therefore, tweaking it may not be allowed. See~\cite{tokcan2025tensor} for an overview of matrix and tensor completion~approaches. \par

Tensor completion problems are typically solved by numerical optimization algorithms (see, e.g.,~\cite{liu2013tc, kressner2013lrtc, song2019tc, pan2021sclrtc, yuan2019ttc, gandy2011tc, bengua2017silrctt}). An~essential aspect of a reliable completion algorithm is its recovery guarantees---specifically, the~conditions under which it can uniquely recover unobserved entries from partial observations. Existing theories generally study probabilistic guarantees for recovery based on conditions such as entries being observed uniformly at random and satisfying incoherence requirements~\cite{candes2010mcvp, candes2012mc, chen2015mc}. However, if~an observation pattern has some structure, better and faster algorithms with deterministic guarantees can be designed by exploiting the structure. Note that, depending on the application, the~observation pattern may be structured rather than random; it may even be fixed, for~instance, when an incomplete dataset is given ``as such'',  without~the possibility of acquiring more entries. There are many interesting problems where conditions like entries being observed uniformly at random may not be explicitly satisfied. In~this work, we discuss one such interesting observation pattern where the fibers of a tensor (along a single specific mode) are either fully observed or entirely missing, unlike the usual entry-wise observations. This observation pattern is interesting because: (i) it occurs in many real-life 
applications; (ii) it makes an intriguing distinction in the uniqueness of the completion between matrix and tensor settings. Specifically, if~some fibers (rows or columns) are entirely missing in a matrix, the~completion problem becomes underdetermined, whereas completion is still possible in a higher-order tensor, even if some fibers are completely missing along a specific mode (see, e.g.,~\cite{chen2018tc, mikael2019fibersamp, stijn2023fibersamp,stijn2023mlsvdfsj}).  In~fact, there are many applications where it is easier to collect data (or sample a multivariate function) along one mode (variable) than along others.  For~example, consider weather time series (e.g., temperature and humidity) collected across various locations~\cite{sofi2024stf}. Think of collecting the data (temperature, latitude, latitude) at specific combinations of geospatial locations. Another example of a fiber-wise observation pattern arises when recording traffic speed data (road segment, day,  time window) for specific combinations of road segments and days (see, e.g.,~\cite{chen2018tc}, where this particular dataset is studied). Other examples include chemical reaction data with ``time'' and ``concentration'' modes. Obtaining samples along the temporal mode may be easier than varying the concentrations (which may require conducting new experiments). \par

Algebraic algorithms exploit the low-rank structure in a specific algebraic manner to design completion algorithms that rely solely on standard numerical linear algebra (NLA) techniques. These algorithms are fast and are guaranteed to work under reasonable deterministic conditions on an observation pattern. In~this line, a~novel algebraic method for fitting a low-rank matrix to a matrix with missing entries was proposed~\cite{jacob2001, hongjun2009mc, stijn2023mlsvdfsj}. More detailed results on low-rank matrix completion, with~extensive discussions on recovery guarantees, have also been presented (see, e.g.,~\cite{kiray2015mc, pimentel2016dsc, stijn2023mlsvdfsj}). Previous studies have shown that the canonical polyadic decomposition (CPD) and multilinear singular value decomposition (MLSVD) of an incomplete tensor, observed along a single mode, can be computed using only standard NLA by exploiting the fiber-wise observation pattern~\cite{mikael2019fibersamp, stijn2023fibersamp, stijn2023mlsvdfsj}. As~big data becomes more prevalent, the~need for stable and scalable algorithms has become more pressing. The~tensor train (TT) decomposition is stable, like the MLSVD and breaks the curse of dimensionality, like the CPD, since it has, asymptotically, the same number of parameters~\cite{oseledets2009breaking, oseledets2010tensortrain}. In~light of these advantages, we present some of the results, including an extension of the algebraic algorithm to the TT format~\cite{shakir2024ttfw}. Note that there is an important difference compared with the popular technique of (TT) cross-approximation, in~the sense that the latter obtains TT decomposition from a subset of fibers sampled across {\em all} modes~\cite{oseledets2010ttcross}. In~our method, fibers are sampled along a {\em single} mode.

\subsection{Contributions}
In this paper, we present more detailed results regarding our approach, with~the following main~additions:
\begin{itemize}
	\item  We provide more insights into piecewise subspace learning, specifically the conditions for determining the column space of a low-rank matrix, where only some pieces (submatrices) are observed.
	\item We include the subspace intersection approach, in~addition to the subspace constraint approach, for~computing the column space of a low-rank matrix from an {\em informationally complete} set of observed submatrices.  
	\item We utilize the piecewise subspace learning approaches to compute TT approximation using only standard NLA operations.
	\item Convincing numerical experiments have been included to show that the proposed method is practically fast and reliable.
	\item In line with~\cite{nico2019compactrep}, we show that the TT decomposition obtained through algebraic completion can serve as a ``proxy'' for efficient subsequent computations (e.g., a~constrained CPD is fitted to the TT approximation rather than to the actual tensor).
\end{itemize}

\subsection{Preliminaries and~Notation} \label{subsec:notation}
We use lower-case, bold lower-case, bold capital, and~calligraphic letters to denote scalars, vectors, matrices, and~tensors, i.e.,~$x$, $\vec{x}, \mat{x},$ and  $\ten{x},$ respectively. The~{\em  order } of a tensor $\ten{x} \in \R^{I_{1} \times I_{2} \times \cdots \times I_{N}}$ is the number of {\em modes} or {\em  ways} it has.  MATLAB-like indexing is used to specify a particular part of a tensor, employing commas, colons, and~semicolons. For~instance, the~specific entry of the $N$th-order tensor is denoted by $x_{i_{1} i_{2} \cdots i_{N}}= \ten{x}(i_{1}, i_{2}, \ldots, i_{N})$ with $i_{n} \in \{ 1, 2, 3, \ldots, I_{n}\}$ for $ n \in \{1, 2, 3, \ldots, N \}$. A~mode-$n$ {\em fiber} of the $N$th-order tensor is obtained by fixing every index except the $n$th. For~a third-order tensor $\ten{x} \in \R^{I_{1} \times I_{2} \times I_{3}}$, the~mode-1 fibers $\vec{x}_{:i_{2} i_{3}}$, mode-2 fibers $\vec{x}_{i_{1} : i_{3}}$, and~mode-3 fibers $\vec{x}_{i_{1} i_{2} :}$ are also known as column, row, and~tube fibers, respectively. Similarly,~mode-1 slices  $\mat{x}_{i_{1} : :}$, mode-2 slices  $\mat{x}_{: i_{2} :}$, and~mode-3 slices  $\mat{x}_{: : i_{3}}$ are also known as horizontal, lateral, and~frontal slices, respectively. The~permutation and reshaping of a tensor are reflected in the ordering of its indices and the position of the semicolon (;), where the semicolon indicates a new mode, and~the ordering of indices indicates the order in which the entries are stacked. In~the $n$th matrix unfolding that is denoted by $\mat{x}_{[1, \ldots, n ; n+1, \ldots, N]} \in \R^{I_1\cdots I_n \times I_{n+1} \cdots I_N}$, the~first $n$ indices enumerate rows and the remaining indices enumerate columns. The~rank, column space (also called the range) and kernel of a matrix $\mat{x}$ are denoted by $\rank{\mat{x}}, \col{\mat{x}}$, and $\ker(\mat{x}),$ respectively. The~dimension of subspace $S$ is denoted by $\dim(S)$. Let $\vec{e}_{i}^{(I)}\in \{0,1\}^{I}$ denote a vector with a unit entry at index $i$ and zeros elsewhere.   We denote a set with a lower-case Greek letter.  The~cardinality of a set $\alpha$ is denoted by $|\alpha|.$ The union and intersection of a sequence of the sets are denoted by  $\cup_{l=1}^L \alpha_{l}$ and  $\cap_{l=1}^L \alpha_{l},$ respectively. \par

Furthermore, we introduce the following notation to represent a submatrix of a matrix $\mat{X} \in \R^{J \times K}$. Let $\alpha_{l} = \{j_1, \ldots, j_{J_l}\} \subseteq \{1, \ldots, J\}$ and $\beta_{l} = \{k_1, \ldots, k_{K_l}\} \subseteq \{1, \ldots, K\}$ be the indices of the $J_l$ selected rows and $K_l$ selected columns, respectively. We define the corresponding row and column selection matrices as $\mat{S}_{r}^{(l)} = \left[
\vec{e}_{j_1}^{(J)}~\ldots ~\vec{e}_{j_{J_{l}}}^{(J)}
\right]^{\top} \in \{0,1\}^{J_l \times J}$ and $\mat{S}_{c}^{(l)} = \left[
\vec{e}_{k_1}^{(K)}~\ldots ~\vec{e}_{k_{K_{l}}}^{(K)}
\right] \in \{0,1\}^{K \times K_l}$, respectively. Matrix $\mat{S}_{r}^{(l)}$ selects the $J_l$ rows indexed by $\alpha_l$ upon left multiplication, and~$\mat{S}_{c}^{(l)}$ selects the $K_l$ columns indexed by $\beta_l$ upon right multiplication. For~example, the~submatrix $\mat{X}_{obs}^{(l)} = \mat{S}_{r}^{(l)} \mat{X} \mat{S}_{c}^{(l)} \in \R^{J_l \times K_l}$ contains elements of $\mat{X}$ observed at $(j, k) \in \alpha_l \times \beta_l$.

\begin{Definition}[Isorank submatrix]
	A submatrix $\mat{X}_{obs}^{(l)}$ of a matrix $\mat{X}$ is called an isorank submatrix or rank-preserving  submatrix if 
	\[ \rank{\mat{X}_{obs}^{(l)}} = \rank{\mat{X}}, \]
	where $\mat{X}_{obs}^{(l)} = \mat{S}_{r}^{(l)} \mat{X} \mat{S}_{c}^{(l)}$. Here, $\mat{S}_{r}^{(l)}$ and $\mat{S}_{c}^{(l)}$ are row and column selection matrices defined by sets $\alpha_l$ and $\beta_l$, respectively, such that $\min(|\alpha_l|, |\beta_l|) \geq \rank{\mat{X}}$.
\end{Definition}

\begin{Definition}[Row overlap]
	Two submatrices $\mat{X}_{obs}^{(1)}$ and $\mat{X}_{obs}^{(2)}$ of a matrix $\mat{X},$ with the corresponding row selection index sets $\alpha_1$ and $\alpha_2$, are said to have a row overlap of size $k$ if there exist $k$ indices that are common to both $\alpha_1$ and $\alpha_2$, i.e.,~$| \alpha_1 \cap \alpha_2| = k$.
\end{Definition}

\begin{Definition}[Contraction] The product between the last mode of tensor $\ten{x} \in \R^{I_1 \times I_2 \times \ldots \times I_N}$ and the first mode of tensor $\ten{y} \in \R^{J_1 \times J_2 \times \ldots \times J_M},$ where $I_N = J_1 = K,$ yields an $(N+M-2)$-order tensor $\ten{z}=\ten{x} \bullet \ten{y},$ the elements of which are given by
	$$
	z_{i_1 \ldots i_{N-1} j_2 \ldots j_M}=\sum_{k=1}^{K} x_{i_1 \ldots i_{N-1} k} y_{k j_2 \ldots j_M}.
	$$ 
\end{Definition}
\subsection{Organization} 
In Section~\ref{sec:tt}, we provide a brief overview of the TT decomposition of a fully observed tensor. Following that, in~Section~\ref{sec:ttfw}, we outline our method for obtaining the TT decomposition of a tensor observed fiber-wise along a single mode. In~Sections~\ref{subsec:algssc} and~\ref{subsec:uniq}, we discuss the algorithm and uniqueness conditions, respectively. Finally, we present convincing numerical experiments in Section~\ref{sec:exp} and a brief conclusion in Section~\ref{sec:con}.

\section{TT Decomposition of Fully Observed~Tensor} \label{sec:tt}%
\begin{figure}[H]
	\centering
	\begin{tikzpicture}[scale=0.6, node distance=0.35cm, chain, all tensors/.append style={tensor scale=0.3, 2d,   no depth compensation}]
		\node[frontal matrix, dim={4,3},  fancy] (G1) {$\mathcal{G}^{(1)}$};
		\node [above  = 0.5cm of G1.center, label= $I_1$]{};
		
		\node [right= 0.35cm of G1, frontal matrix, dim={1,0.5},  fill=kulblue70] (r1) {};
		\node [above  = 0.1cm of r1.center, left = 0.1cm of r1.center ,  label = $R_{1}$]{};
		
		\node [right=0.35cm of r1,  tensor, dim={2,5,3}, fancy] (G2){$\mathcal{G}^{(2)}$};
		\node [above  =  -0.1cm of G2.top face center, label= $I_2$]{};
		
		\node [right=0.35cm of G2, frontal matrix, dim={1,0.5},  fill=kulblue70] (r2) {};
		\node [above  = 0.1cm of r2.center, left = 0.1cm of r2.center ,  label = $R_{2}$]{};
		
		\node [right=0.35cm of r2,tensor, dim={3,6,3},  fancy](G3){$\mathcal{G}^{(3)}$};
		\node [above  =  -0.1cm of G3.top face center, label= $I_3$]{};
		
		\node [right=0.35cm of G3, frontal matrix, dim={1,0.5},  fill=kulblue70] (r3) {};
		\node [above  = 0.1cm of r3.center, left = 0.1cm of r3.center ,  label = $R_{3}$]{};
		
		\node [right=0.35cm of r3,tensor, dim={3,4,2},  fancy] (G4){$\mathcal{G}^{(4)}$};
		\node [above  =  -0.1cm of G4.top face center, label= $I_4$]{};
		
		\node [right=0.35cm of G4, frontal matrix, dim={1,0.5}, fill=kulblue70] (r4) {};
		\node [above  = 0.1cm of r4.center, left = 0.1cm of r4.center ,  label = $R_{4}$]{};
		
		\node[right=0.35cm of r4,frontal matrix, dim={2,4},  fancy] (G5) {$\mathcal{G}^{(5)}$};
		\node [above  = 0.2cm of G5.center, label= $I_5$]{};
		\draw [-] (G1.east) -- (r1.west);
		\draw[-](r1.east) -- (G2.west);
		\draw[-](G2.east) -- (r2.west);
		\draw[-](r2.east) -- (G3.west);
		\draw[-](G3.east) -- (r3.west);
		\draw[-](r3.east) -- (G4.west);
		\draw[-](G4.east) -- (r4.west);
		\draw[-](r4.east) -- (G5.west);
	\end{tikzpicture}
	\caption{TT decomposition of $5$th-order tensor as a train of five core tensors, where $\ten{g}^{(1)} \in \R^{1 \times I_{1} \times R_{1}}$ and $\ten{g}^{(5)} \in \R^{R_{4} \times I_{5} \times 1}.$} 
	\label{fig:tt}
\end{figure}
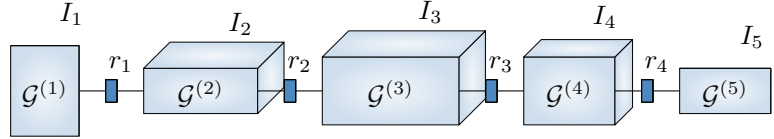

A TT decomposition of a tensor $\ten{x} \in \R^{I_1 \times I_2 \times \cdots \times I_N}$ corresponds to a contraction of a sequence of third-order core tensors (TT cores) $\ten{g}^{(n)} \in \R^{R_{n-1} \times I_n  \times R_{n}} (1 \leq n \leq N)$ with $R_{0} = R_{N} =1$, such that each entry of $\ten{x}$ can be expressed as the sequence of matrix products~\cite{verstraete2007MPS, oseledets2010tensortrain}:
\begin{equation}
	\label{eqn:tt1}
	\ten{x} = \ten{g}^{(1)} \bullet \ten{g}^{(2)} \bullet \cdots \bullet \ten{g}^{(N-1)} \bullet \ten{g}^{(N)},
\end{equation}
or entry-wise, we can write:
\begin{equation}
	\label{eqn:tt2}
	x_{i_{1} i_{2} \cdots i_{N-1} i_{N}} = \mat{g}^{(1)}_{: i_{1} :} \mat{g}^{(2)}_{: i_{2} :}  \cdots \mat{g}^{(N-1)}_{: i_{N-1} :} \mat{g}^{(N)}_{: i_{N} :},
\end{equation}
where the matrix $\mat{g}^{(n)}_{: i :} \in \R^{R_{n-1} \times R_{n}}$ is the $i$th mode-2 slice of the core tensor $\ten{g}^{(n)}.$ The tuple of minimal integers $(R_0, \ldots, R_{N})$ for which the equality in Equations \eqref{eqn:tt1} and \eqref{eqn:tt2} holds is the TT rank of $\ten{x},$ denoted by $\text{rank}_{TT}(\ten{x}).$ The dense tensor $\ten{x}$ has a total storage complexity of $\prod_{n = 1}^{N} I_{n}$, whereas in the TT format, the~storage complexity is $\sum_{n = 1}^{N} R_{n-1}I_{n}R_{n}$. Hence, a~low TT rank can greatly reduce storage complexity. A~visualization of the TT decomposition of a $5$th-order tensor is shown in Figure~\ref{fig:tt}.

Additionally, we introduce a shorthand notation for TT decomposition using partial products. A~tensor $\ten{x}$ can be represented as
\begin{equation}
	\label{eqn:tt3}
	\ten{x} = \ten{g}^{(<n)} \bullet \ten{g}^{(n:m)} \bullet \ten{g}^{(>m)}, \quad 1 \leq n \leq  m \leq N, 
\end{equation}
where the left partial product $\ten{g}^{(<n)} = \ten{g}^{(1)} \bullet  \cdots \bullet \ten{g}^{(n-1)},$ the right partial product $\ten{g}^{(>m)} = \ten{g}^{(m+1)}  \bullet \cdots \bullet \ten{g}^{(N)}$ and $\ten{g}^{(n:m)} = \ten{g}^{(n)} \bullet \cdots \bullet \ten{g}^{(m)}.$ \par

TT decomposition of a given tensor $\ten{x}$ is computed by a sequence of truncated SVDs (TT-SVD)~\cite{oseledets2010tensortrain}. The~TT-SVD algorithm sequentially computes TT cores and essentially alternates between the SVD computation step and projection step, which makes it difficult to parallelize. Inspired by the natural parallelizability of MLSVD~\cite{delathauwer2000mlsvd}, a~similar method was proposed that computes the orthonormal basis for the column spaces of all matrix unfoldings independently and then utilizes them to compute the TT cores (Parallel-TTSVD)~\cite{shi2023paralleltt}. It should be noted that independently computing the orthonormal basis for the column spaces of each matrix unfolding is more expensive than the basis computations in \mbox{TT-SVD,} where the projection step reduces the complexity of subsequent computations.  One can reduce this computational burden by replacing SVDs with randomized SVDs~\cite{mahoney2010randomizedalgorithms, halko2011randsvd}. The~algorithm for computing the TT decomposition with parallel SVDs (Parallel-TTSVD) is shown in Algorithm~\ref{alg:ptt} (for more in-depth details, refer to \cite{shi2023paralleltt}, Algorithm 3.1). 

\vspace{6pt}
\begin{algorithm}[H]
	\caption{Parallel-TTSVD}\label{alg:ptt}
	\Indm
	\KwIn{A tensor $\ten{x} \in \R^{I_{1} \times I_{2} \times \cdots \times I_{N}},$ and TT rank $(R_{0}, \ldots , R_{N})$\;}
	\KwOut{TT cores $\ten{g}^{(1)}, \ldots, \ten{g}^{(N)}$\;}   
	
	\Indp   
	\tcc{compute orthonormal bases for the ranges of the $n$th unfoldings}
	\For{$ 1 \leq n \leq N-1$}{
		Compute $\mat{U}^{(n)}, \mat{\Sigma}^{(n)}, \mat{V}^{(n) \top} = \texttt{svd}\left(\mat{X}_{[1, \ldots,  n ; n+1, \ldots, N]},  R_{n}\right)$\; \label{line:rangefull}
		Set $\mat{A}^{(n)} = \mat{U}^{(n)}, $ and $\mat{B}^{(n)} = \mat{\Sigma}^{(n)}\mat{V}^{(n)\top}$ \;
	}
	\tcc{compute the TT cores}
	Set $\ten{g}^{(1)} = \texttt{reshape}\left(\mat{a}^{(1)}, [1, I_1, R_1]\right)$ \;
	\For{$ 1 \leq n \leq N-2$}{ 
		Compute $\mat{w}^{(n)} = \mat{a}^{(n)\top} \left(\mat{a}^{(n+1)}\right)_{[1, \ldots, n ; n+1, (n+1)+1]}$\;
		Set $\ten{g}^{(n+1)} = \texttt{reshape}\left(\mat{w}^{(n)}, [R_{n}, I_{n+1}, R_{n+1}] \right)$\;
	}
	Set $\ten{g}^{(N)} = \texttt{reshape}\left(\mat{b}^{(N-1) \top}, [R_{N-1}, I_N, 1]\right)$ \;
\end{algorithm}

\section{TT Completion from Fiberwise Observations Along a Single Mode} \label{sec:ttfw}
{\textbf{Roadmap:} This section presents our main contribution. Our goal is to design an algebraic algorithm for computing the TT decomposition---analogous to Algorithm~\ref{alg:ptt}, but~able to handle an incomplete tensor observed through mode-$N$ fibers. We begin by examining the structured observation pattern that appears in the matrix unfoldings of such tensors. In~Section~\ref{subsec:pssl}, we show how to compute the column space of a low‑rank matrix formed by stacking two slices of a fiber-wise observed tensor, and~then extend this approach to the partially observed matrix unfolding that we are actually dealing with (and which consists of multiple slices). This leads to two subspace‑estimation methods presented in Sections~\ref{subsec:ssc} and~\ref{subsec:ssI}.  Finally, in~Section~\ref{subsec:algssc}, we use this core subspace‑estimation step to build the full TT completion algorithm and, in~Section~\ref{subsec:uniq}, summarize the conditions under which it is guaranteed to succeed.}

\begin{figure}[H]
	\centering
	\begin{tikzpicture}[scale=2]
		\begin{scope}[xshift=0cm,yshift=0cm]
			\node at (4.9, 9.5){$\mat{x}_{[1, \ldots,  n ; n+1, \ldots, N]} ~= $};
		\end{scope}
		
		\begin{scope}[xshift=1cm,yshift=0cm]
			\node at (6.5, 9.5){\rotatebox{0}{\scalebox{1.5}{$\cdots$}}};
			
			\draw [ fill = white] (5,10) rectangle (6,9);
			\draw [ ] (6,10) rectangle (7,9);
			\draw [ fill = white] (7,10) rectangle (8,9);
			
			\draw [ fill = kulblue90] (5,9.9) rectangle (6,9.8);
			\draw [ fill = kulblue90] (5,9.5) rectangle (6,9.4);
			\draw [ fill = kulblue90] (5,9.2) rectangle (6,9.1);
			
			\draw [ fill = kulblue90] (7,9.6) rectangle (8,9.5);
			\draw [ fill = kulblue90] (7,9.4) rectangle (8,9.3);
			\draw [ fill = kulblue90] (7,9.1) rectangle (8,9);
		\end{scope}
	\end{tikzpicture}
	\caption{\small{The $n$th matrix unfolding of a tensor observed through fibers along a single mode is characterized by a structured observation pattern.}}
	\label{fig:slicestacks}
\end{figure}
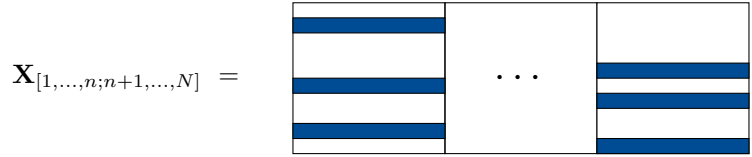

Assume a given tensor $\ten{x} \in \R^{I_{1} \times \cdots \times I_{N}}$ with $\text{rank}_{TT}(\ten{x}) = (R_{0}, \ldots, R_{N})$ of which only some of the mode-$N$ (last mode) fibers are observed. One essential step to compute the TT cores is to find orthonormal bases for the ranges of the matrix unfoldings.  However, since the matrix unfoldings have missing entries, it is not directly possible to utilize, e.g.,~SVD or QR factorization to obtain the ranges.  Note that under fiber-wise observations, two types of observation patterns arise in matrix unfoldings. In~the $(N-1)$th matrix unfolding, i.e.,~$\mat{x}_{[1, \ldots, N-1 ; N]}$, the~rows are either fully observed or entirely missing. If~the number of observed rows is greater than or equal to $R_{N-1}, $ we can generically (i.e., with probability 1 when the matrix entries are drawn from a continuous distribution) obtain the last TT core $(\ten{g}^{(N)} \in \R^{R_{N-1} \times I_N  \times 1})$  by computing the top $R_{N-1}$ right-singular vectors of the matrix formed by the observed rows. The~other matrix unfoldings can be seen as the horizontal stack of mode-2 slices $\tilde{\mat{x}}_{:i:} $ of a third-order reshaping $\tilde{\ten{x}} =\ten{x}_{[1, \ldots, n ; n+1, \ldots, N-1 ; N]}\in \R^{\prod_{i = 1}^{n} I_{i} \times \prod_{i = n+1}^{N-1} I_{i} \times I_{N}} $ of the tensor $\ten{x}.$ That is,
\begin{equation}
	\label{eqn:ttfw1}
	\mat{x}_{[1, \ldots,  n ; n+1, \ldots, N]} = \left[ \tilde{\mat{x}}_{:1:} ~\cdots ~ \tilde{\mat{x}}_{:\mathsmaller{\prod_{i = n+1}^{N-1} I_{i}}:} \right]. 
\end{equation}
In this case, the~rows of the lateral slices (submatrices of the unfolding) of  $\tilde{\ten{x}}$ are either fully observed or entirely missing. However, different slices may have observed rows at different row indices. A~visual representation of such a matrix unfolding is shown in  Figure~\ref{fig:slicestacks}.

The following section will discuss how to determine the orthonormal basis for the column space of such a partially observed matrix~unfolding.

\subsection{Piecewise Subspace~Learning}\label{subsec:pssl}
In this section, we will discuss how to determine the overall column space of a low-rank matrix, of~which only some pieces (submatrices) are observed. Moreover, we will characterize the sampling of pieces such that the overall subspace is guaranteed to be unique (i.e., the~subspace is identifiable). Subspace identifiability is closely related to the matrix completion problem. The~conditions required to identify the subspace (e.g., column space) are necessary but not sufficient for matrix completion~\cite{pimentel2015dsc, stijn2023mlsvdfsj}. We will first explore the conditions for unique matrix completion, a~problem that is widely studied in the literature, and~then discuss how these conditions are relaxed in the context of subspace~identification.  \par
For simplicity, we use $\mat{m} \in \R^{J \times K}$ to denote a rank-$R$ matrix unfolding of which we need to determine the $R$-dimensional column space. In~the literature, unique minimal rank completions for various types of observation patterns have been studied (see, e.g.,~ \cite{bostian2001uniq, kaashoek1988uniq, epperly2021minrank}). This work will focus on an observation pattern shown in Figure~\ref{fig:slicestacks}.  Before~proceeding to the general case, we will begin with a simple  matrix formed by stacking two slices of a tensor, as~explained next.  Denote the observed submatrices of $\mat{m}$ by $\mat{m}_{obs}^{(1)} \in R^{(J_1 + J_2) \times K_1}$  and $\mat{m}_{obs}^{(2)} \in R^{(J_2 + J_3) \times K_2}$ that have $J_2$ overlapping rows, and~$J=J_1+J_2+J_3$ and $K=K_1 + K_2$. (Observed submatrices of $\mat{M}$ are shown shaded.)   While the overlap can occur anywhere among the rows, let us assume, without~loss of generality,  that the row overlap occurs in the middle block, as~ shown below: 
\begin{equation}
	\label{eqn:mconn}
	\mat{M} = 
	\left[ \begin{array}{c|c}
		\cellcolor{kulblue20} \mat{M}_{11} & \mat{M}_{12} \\
		\cellcolor{kulblue20} \mat{M}_{21} & \cellcolor{kulblue20} \mat{M}_{22} \\
		\mat{M}_{33} & \cellcolor{kulblue20} \mat{M}_{32}
	\end{array} \right], \quad 
	\mat{m}_{obs}^{(1)} = \left[ \begin{array}{c}
		\mat{M}_{11} \\
		\mat{M}_{21}
	\end{array} \right] \text{ and }  \mat{m}_{obs}^{(2)} = \left[ \begin{array}{c}
		\mat{M}_{22} \\
		\mat{M}_{32}
	\end{array} \right]. 
\end{equation}
It has been proven that there is a unique rank-$R$ completion for $\mat{M}$ if and only if the following condition holds (see, e.g.,~Corollary 2.3 in~\cite{bostian2001uniq}):
\begin{equation}
	\label{cond:overlap}
	\rank{\mat{m}_{obs}^{(1)}} = \rank{\mat{M}_{21}} =
	\rank{\left[\mat{M}_{21} ~ \mat{M}_{22} \right]}  = \rank{\mat{M}_{22}} = \rank{\mat{m}_{obs}^{(2)}} = R.
\end{equation}

As noted above, subspace identifiability is less restrictive than matrix completion. To~illustrate this, consider an incomplete matrix formed by concatenating a matrix $\mat{M}$ that satisfies condition~\eqref{cond:overlap}, and~a vector $\vec{z}$. Let this concatenated matrix be denoted as \mbox{$\hat{\mat{M}} = \left[\mat{M} \mid \vec{z}\right]$}, and~assume that $\rank{\hat{\mat{M}}} = \rank{\mat{M}} = R$. The~subspace identifiability of  $\hat{\mat{M}}$ is ensured by $\mat{M}$ since the latter satisfies the row overlapping condition. As~a result, $\col{\hat{\mat{M}}}$ can be determined solely from $\mat{M}$~\cite{pimentel2015dsc, pimentel2016dsc, stijn2023mlsvdfsj}. However, ensuring the unique completion of $\hat{\mat{M}}$ requires an additional condition: the vector $\vec{z}$ should also have at least $R$ observed entries. Without~this, a~unique recovery is not~possible. \par

The study in~\cite{stijn2023mlsvdfsj} investigates the subspace identifiability of a partially observed low-rank matrix formed by stacking multiple slices (submatrices) of an incomplete tensor as shown in Figure~\ref{fig:slicestacks}. The~authors show that subspace identifiability also requires a row-overlapping condition as stated in~\eqref{cond:overlap}. However, this is not very restrictive as only some, rather than all, pairs of observed submatrices need to satisfy it. As~noted earlier, for~the full completion to be unique, we also need to ensure that every column contains at least $R$ observed entries. For~further details, see~\cite{pimentel2015dsc, pimentel2016dsc, stijn2023mlsvdfsj}.

In the following example of a partially observed rank-$1$ matrix, we provide a geometrical interpretation of the conditions required to identify the column space and set out the basis for the algebraic algorithm for computing the desired column space more~generally.  

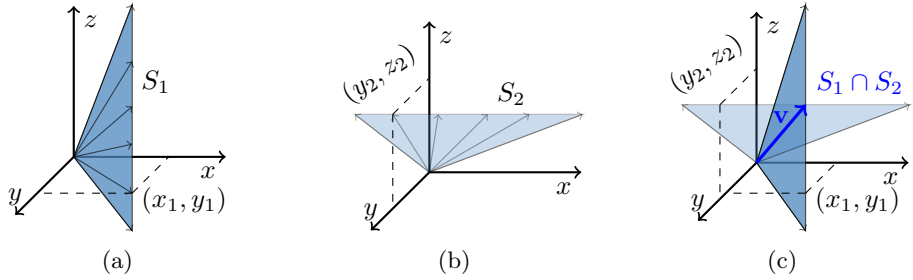
\begin{figure}[H]
	\begin{minipage}[b]{0.3\textwidth}
		\centering
		\begin{tikzpicture}
			\draw[thick, ->] (0,0,0) -- (2,0,0) node[anchor=north east]{$x$};
			\draw[thick, ->] (0,0,0) -- (0,2,0) node[anchor=north west]{$z$};
			\draw[thick, ->] (0,0,0) -- (0,0,2) node[anchor=south]{$y$};
			
			\draw[fill=kulblue70,opacity=0.75] (1.25,-0.5,1.25) -- (1.25,0,1.25) -- (1.25, 0.65,1.25) -- (1.25, 1.15,1.25)--(1.25, 1.75,1.25) -- (1.25,2.5,1.25) -- (0,0,0) -- cycle;
			
			\draw[dashed] (1.25,0,1.25) -- (1.25,0,0);
			\draw[dashed] (1.25,0,1.25) -- (0,0,1.25);
			
			\node[anchor=north west] at (1.25,0.2,1.25) {$(x_1, y_1)$};
			\draw[->, fill=black,opacity=0.65] (0,0,0) -- (1.25, -0.5,1.25);
			\draw[->, fill=black,opacity=0.65] (0,0,0) -- (1.25, 0,1.25);
			\draw[->, fill=black,opacity=0.65] (0,0,0) -- (1.25, 0.65,1.25);
			\draw[->, fill=black,opacity=0.65] (0,0,0) -- (1.25, 1.15,1.25);
			\draw[->, fill=black,opacity=0.65] (0,0,0) -- (1.25, 1.75,1.25) node[anchor=north west, text=black, opacity=1]{$S_1$};
			\draw[->, fill=black,opacity=0.65] (0,0,0) -- (1.25, 2.5,1.25);
		\end{tikzpicture}
	\end{minipage}
	\begin{minipage}[b]{0.3\textwidth}
		\centering
		\centering
		\begin{tikzpicture}
			\draw[thick, ->] (0,0,0) -- (2,0,0) node[anchor=north east]{$x$};
			\draw[thick, ->] (0,0,0) -- (0,2,0) node[anchor=north west]{$z$};
			\draw[thick, ->] (0,0,0) -- (0,0,2) node[anchor=south]{$y$};
			
			\draw[fill=kulblue70,opacity=0.3] (-0.5,1.25,1.25) -- (0,1.25,1.25) -- (0.6,1.25,1.25) -- (1.25,1.25,1.25) -- (1.8,1.25,1.25) -- (2.5,1.25,1.25) -- (0,0,0) -- cycle;
			
			\draw[dashed] (0,1.25,1.25) -- (0,1.25,0);
			\draw[dashed] (0,1.25,1.25) -- (0,0,1.25);
			
			\node[anchor=east, rotate = 35] at (0,1.8,0.3) {$(y_2, z_2)$};
			\draw[->, fill=black,opacity=0.35] (0,0,0) -- (-0.5,1.25,1.25);
			\draw[->, fill=black,opacity=0.35] (0,0,0) -- (0,1.25,1.25);
			\draw[->, fill=black,opacity=0.35] (0,0,0) -- (0.6,1.25,1.25);
			\draw[->, fill=black,opacity=0.35] (0,0,0) -- (1.25,1.25,1.25) node[anchor=south west, text=black, opacity=1]{$S_2$};
			\draw[->, fill=black,opacity=0.35] (0,0,0) -- (1.8,1.25,1.25);
			\draw[->, fill=black,opacity=0.35] (0,0,0) -- (2.5,1.25,1.25);
		\end{tikzpicture}
	\end{minipage}
	\begin{minipage}[b]{0.3\textwidth}
		\centering
		\begin{tikzpicture}
			\draw[thick, ->] (0,0,0) -- (2,0,0) node[anchor=north east]{$x$};
			\draw[thick, ->] (0,0,0) -- (0,2,0) node[anchor=north west]{$z$};
			\draw[thick, ->] (0,0,0) -- (0,0,2) node[anchor=south]{$y$};
			
			\draw[fill=kulblue70,opacity=0.3] (-0.5,1.25,1.25) -- (2.5,1.25,1.25) -- (0,0,0) -- cycle;
			\draw[fill=kulblue70,opacity=0.75] (1.05,-0.5,1.05) -- (1.05,2.5,1.05) -- (0,0,0) -- cycle;
			
			\draw[dashed] (0,1.25,1.25) -- (0,1.25,0);
			\draw[dashed] (0,1.25,1.25) -- (0,0,1.25);
			\node[anchor=east, rotate = 35] at (0,1.8,0.3) {$(y_2, z_2)$};
			\draw[dashed] (1.05,0,1.05) -- (1.05,0,0);
			\draw[dashed] (1.05,0,1.05) -- (0,0,1.05);
			\node[anchor=north west] at (1.05,0.2,1.05) {$(x_1, y_1)$};
			
			\draw[->, fill=black,opacity=0.68] (0,0,0) -- (1.05, -0.5,1.05);
			\draw[->, fill=black,opacity=0.68] (0,0,0) -- (1.05, 2.5,1.05);
			
			\draw[->, fill=black,opacity=0.35] (0,0,0) -- (-0.5,1.25,1.25);
			\draw[->, fill=black,opacity=0.35] (0,0,0) -- (2.5,1.25,1.25);
			
			\draw[->, very thick, blue] (0,0,0) -- (1.13,1.25,1.25) node[pos=0.5, above]{$\vec{v}$} node[pos=1, above right]{$ S_1 \cap S_2$};
			
		\end{tikzpicture}
	\end{minipage}\\
	\centering (\textbf{a})\hspace{110pt}(\textbf{b})\hspace{110pt}(\textbf{c}) 
	\caption{(\textbf{a}) Affine subspace $S_1$ corresponding to the first column $\vec{m}_{:1}$; (\textbf{b}) affine subspace $S_2$ corresponding to the second column $\vec{m}_{:2}$; (\textbf{c}) intersection of the affine subspaces $S_1 \cap S_2$.}
	\label{fig:affine}
\end{figure}

\begin{Example} \label{ex:1} A partially observed rank-$1$ matrix is given: $\mat{M} = \left[ \begin{array}{c | c}
		\cellcolor{kulblue50} x_1 & x_2 \\
		\cellcolor{kulblue50} y_1 & \cellcolor{kulblue20} y_2 \\
		z_1 & \cellcolor{kulblue20} z_2
	\end{array} \right].$ Suppose we are required to determine its 1-dimensional range in $\mathbb{R}^{3}.$
	
	We discuss the subspace identifiability of the rank-1 matrix $\mat{m}$ in terms of subspaces associated with partially observed columns. In~the first column, $z_1$ is missing and may lay anywhere on the line parallel to the z-axis passing through $(x_1, y_1),$ as shown in Figure~\ref{fig:affine}a. The~subspace $S_1$ (shaded region) accounts for all possible completions of the first column. Similarly, in~the second column,~entry $x_2$  is missing and may lay anywhere on the line parallel to the x-axis passing through $(y_2, z_2)$, and the corresponding subspace that accounts for all possible completions of the second column is  $S_2$ (shaded region), as~shown in Figure~\ref{fig:affine}b. Given that the unknown 1-dimensional subspace $\col{\mat{M}}$ spans both columns, it must lie in the intersection of these affine subspaces, i.e.,~$\col{\mat{M}} \subseteq S_1 \cap S_2$. Note that there is a subset relationship between $\col{\mat{m}}$ and $S_1 \cap S_2$, not an equality. However, if~\( \dim(S_1 \cap S_2) = 1 \), i.e.,~if \( S_1 \cap S_2 \) forms a line, then the desired \( \col{\mat{m}} \) is necessarily equal to \( S_1 \cap S_2 \). Therefore, any nonzero vector \( \vec{v} \in S_1 \cap S_2 \) spans \( \col{\mat{m}} \), as~shown in Figure~\ref{fig:affine}c.
	In this case, any technique that computes the intersection of subspaces can be employed to find the desired subspace $\col{\mat{m}}.$  \par
	One approach to determining the intersecting subspace is via the null spaces. Specifically, find $\vec{n_1}$ and $\vec{n_2}$ such that
	$ \langle \vec{n_1}, \mathbf{m}_{:1} \rangle=0, $  and $ \langle \vec{n_2}, \vec{m}_{:2} \rangle=0.$
	Define $\mat{N} = [\vec{n_1} ~ \vec{n_2}] \in \mathbb{R}^{3 \times 2}$. If~$\dim(\ker{(\mat{N}^{\top})})=1$, then the subspace $S$ is given by $\ker{(\mat{N}^{\top})}$. The~ condition $\dim(\ker{(\mathbf{N}^{\top})})=1$ represents the dual of the condition $\dim(S_1 \cap S_2)=1$. Several other approaches have been proposed to find the intersection of subspaces (see, e.g.,~\cite{benisrael2013proj}). 
\end{Example} 

\vspace{-12pt}

Before proceeding to the rank-$R$ case, for~ease of notation, we use $\mat{m} \in \R^{J \times K}$ to denote a low-rank matrix $\mat{x}_{[1, \ldots,  n ; n+1, \ldots, N]}$ for which we need to determine the $R$-dimensional column space. Let $\mat{m}_{obs}^{(l)} = \mat{S}_{r}^{(l)} \mat{m} \mat{S}_{c}^{(l)} \in \R^{J_l \times K_l}$ represent the $l$th fully observed isorank submatrix of $\mat{m}$. We denote a matrix that holds an orthonormal basis for $\col{\mat{m}}$ by $\mat{a} \in \R^{J \times R}$. Below, we discuss two approaches to determine $\mat{A}$: (1) the {\em subspace constraint} approach, which computes it via null spaces of observed submatrices, and~(2) the {\em subspace intersection} method, which computes it via column spaces of observed~submatrices. \par

\subsubsection{Subspace Constraint~Approach}\label{subsec:ssc}
The main idea of the subspace constraint method is to determine the overall range of a matrix from a set of constraints derived from the submatrices~\cite{jacob2001}. This method involves identifying an informationally complete set of fully observed isorank submatrices $\{\mat{m}_{obs}^{(l)}\}_{l=1}^{L}$ of $\mat{M}$ \cite{stijn2023fibersamp}. By~{\em informationally complete,} we mean that the set of constraints derived from the submatrices should fully characterize the desired subspace $\col{\mat{A}}$. As~we assume $\mat{m}_{obs}^{(l)}$ is an isorank submatrix of $\mat{m},$ we can write $\mat{m}^{(l)}_{obs} = \mat{a}^{(l)}_{obs} \mat{b}^{(l)\top}_{obs}, $ where matrices $\mat{a}^{(l)}_{obs} \in \R^{J_l \times R}$ and  $\mat{b}^{(l)}_{obs} \in \R^{K_l \times R}$  are full column rank. Stacking a basis for the orthogonal complement of $\col{\mat{m}^{(l)}_{obs}}$ in a matrix $\mat{n}_{obs}^{(l)} \in \R^{J_l \times (J_l - R)},$ each column of the resulting matrix after zero padding, i.e.,~ $\mat{s}_r^{(l)\top}\mat{n}_{obs}^{(l)} \in \R^{J \times (J_l - R)}, $ represents a vector that is orthogonal to the desired subspace $\col{\mat{A}},$ thus imposing constraints on what $\col{\mat{A}}$ may be. To~obtain the desired $\col{\mat{A}}$, we find orthogonal complements for $L$ such submatrices $\mat{m}_{obs}^{(l)}$ and concatenate them in a matrix  $\mat{N} = \left[ \mat{s}_r^{(1)\top}\mat{n}_{obs}^{(1)}, \cdots,  \mat{s}_r^{(L)\top}\mat{n}_{obs}^{(L)}\right] \in \R^{J \times (\sum_{l=1}^{L}(J_l-R))}.$ If the $L$ orthogonal complements are enough to ensure the $\ker(\mat{n}^{\top})$ is of minimal dimension, i.e.,~$\dim( \ker(\mat{n}^{\top})) = R,$  where the very existence of $\col{\mat{A}}$ implies that the dimension is at least $R,$ then the desired $\col{\mat{A}}$ is necessarily given by $\ker(\mat{n}^{\top})$.  In~other words, the~number of independent constraints must be at least \( J - R \), which leads to the inequality \( \sum_{l=1}^L (J_l - R) \geq J - R.\) For example, if~we assume that each observed submatrix imposes one (independent) constraint and we have \( J_l = R + 1 \), then the minimum number of the observed submatrices required is \( L = J - R \). This can be seen as the generalization of the rank-$1$ case mentioned in Example~\ref{ex:1}. In~what follows, we assume that our working assumptions hold, i.e.,~$\{\mat{m}_{obs}^{(l)}\}_{l=1}^{L}$ is a set of isorank submatrices of $\mat{m}$ and $\dim(\ker(\mat{n}^{\top})) = R$. Then, this set of submatrices is {\em informationally complete}. In~the noisy case, we estimate $\mat{n}^{(l)}_{obs}$ from the left-singular vectors corresponding to the smallest singular values of $\mat{m}^{(l)}_{obs}$. Let the SVD of $\mat{m}^{(l)}_{obs}$ be $\mat{m}^{(l)}_{obs} = \mat{u}^{(l)} \mat{\Sigma}^{(l)} \mat{v}^{(l) \top}$. Then, $\mat{n}^{(l)}_{obs} = \left[\vec{u}_{:R+1}^{(l)} \cdots \vec{u}_{:J_{l}}^{(l)} \right].$ In the noisy case,  we also estimate the desired $\col{\mat{a}}$ to be the subspace that is most orthogonal to the space spanned by the columns of $\mat{n},$ i.e.,~we estimate $\col{\mat{a}}$ as the (approximate) kernel of $\mat{n}^{\top}.$ Let the SVD of $\mat{n}$ be $\mat{u} \mat{\Sigma} \mat{v}^{ \top}$. Then $\mat{A} = \left[\vec{u}_{:J-R+1} \cdots \vec{u}_{:J} \right].$

\subsubsection{Subspace Intersection~Approach}\label{subsec:ssI}%
Let us first define a binary matrix $\tilde{\mat{S}}_r^{(l)} \in \{0,1\}^{(J-J_l) \times J}$, in~which the rows correspond to standard unit vectors that indicate which rows of $\mat{m}$ are missing in $\mat{m}_{obs}^{(l)}$. (The matrix  $\tilde{\mat{S}}_r^{(l)}$ is such that $\mat{I}_J=\left[ \mat{S}_r^{(l)\top} \quad \tilde{\mat{S}}_r^{(l)\top}\right] \mat{\Pi} \in \R^{J \times J}$ for some permutation matrix $\mat{\Pi} \in \R^{J \times J}$.) 
Additionally, let $S_l$ denote a subspace that accounts for all possible completions of rows that are missing in $\mat{m}_{obs}^{(l)}$.  Specifically, let the SVD of $\mat{m}^{(l)}_{obs}$ be $\mat{m}^{(l)}_{obs} = \mat{u}^{(l)} \mat{\Sigma}^{(l)} \mat{v}^{(l) \top}.$  Then, the~subspace $S_l$ is given by
$ S_l = \operatorname{span}\left( \mat{q}_{S_l} \right),$ where $\mat{q}_{S_l}  \mathrel{\mathop:} = \left[ \mat{S}_r^{(l)\top} [\vec{u}_{:1}^{(l)} \cdots \vec{u}_{:R}^{(l)}] \quad \tilde{\mat{S}}_r^{(l)\top} \right] \in \R^{J \times (R+J-J_l)}$; note that the columns of $\mat{q}_{S_l}$ form an orthonormal basis for the subspace $S_l$.

As illustrated in Example~\ref{ex:1}, the~desired subspace $\col{\mat{M}}$ satisfies $\col{\mat{M}} \subseteq \left(\cap_{l=1}^L S_l \right)$. However, since we assume the subspaces are associated with an informationally complete set of fully observed isorank submatrices $\{\mat{m}_{obs}^{(l)}\}_{l=1}^{L}$, i.e.,~$\dim{\left(\cap_{l=1}^L S_l\right)} = R \Leftrightarrow \dim(\ker(\mathbf{N}^{\top})) = R$, it follows that the $R$-dimensional subspace $\col{\mat{M}}$ is necessarily equal to $\cap_{l=1}^L S_l$. Hence, we simply need to find the intersection of the~subspaces. \par 

A closed-form solution for computing the intersection of $L\geq 2$ subspaces in finite-dimensional spaces is discussed in~\cite{benisrael2013proj}, and~a low-complexity implementation of these formulas based on SVD is presented in~\cite{fenggang2019projK, stijn2023mlsvdfsj}.  The~procedure used to compute the intersection is as follows. Given a set of matrices $\left\{ \mat{Q}_{S_l} \right\}_{l=1}^L,$ the intersection of the corresponding subspaces can be computed as $\cap_{l=1}^L S_l = \ker \left( L\mat{I}_J - \sum_{l=1}^L \mat{Q}_{S_l}\mat{Q}^{\top}_{S_l} \right)$. It has been shown that if $\dim{\left(\cap_{l=1}^L S_l\right)} = R,$ then the solution can be computed more efficiently via the SVD of a concatenated matrix $\mat{Q} = \left [ \mat{Q}_{S_1}, \cdots, \mat{Q}_{S_L} \right] \in \R^{J \times \left(LR + \sum_{l=1}^L(J-J_l)\right)},$ without first calculating the orthogonal projectors $\left\{ \mat{Q}_{S_l}\mat{Q}^{\top}_{S_l} \right\}$. Specifically, let the SVD of \mbox{$\mat{Q}$ be $\mat{Q} = \mat{U \Sigma V}^{\top}$}. Then, $\mat{A} = \left[\vec{u}_{:1}, \cdots, \vec{u}_{:R} \right]$ forms an orthonormal basis for the desired subspace $\col{\mat{M}}$. See~\mbox{\cite{benisrael2013proj, sorensen2021gcca, stijn2023mlsvdfsj}} for an in-depth discussion. In~case the matrix $\mathbf{Q}$ is too large to fit in memory, one can use an incremental SVD to find the dominant left singular subspace~\cite{brand2002iSVD, brand2006iSVD}.

\subsection{Algorithm}\label{subsec:algssc}

This section presents our algorithm for computing the TT decomposition of a tensor $\ten{x} \in \R^{I_1 \times I_2 \times \cdots \times I_N}$ observed along the $N$th mode. Our algorithm is similar to Algorithm~\ref{alg:ptt}, with~the key subspace computation step in line~\ref{line:rangefull} adapted to handle a tensor observed through mode-$N$ fibers. The~key steps of the algorithm are as follows. The~piecewise subspace learning approach is employed to compute orthonormal bases for the column spaces of the partially observed matrix unfoldings \(\mat{X}_{[1,\ldots,n; \, n+1,\ldots, N]} \) for \(1 \leq n \leq N\!-\!2\), assuming that the corresponding observed submatrices satisfy the informational-completeness condition of Section~\ref{subsec:uniq}. These orthonormal bases are used to compute the TT cores \(\ten{G}^{(1)}, \ldots, \ten{G}^{(N\!-\!2)}\). The~last core $\ten{g}^{(N)}$ is obtained from an orthonormal basis of the observed rows of the $(N\!-\!1)$th matrix unfolding, and~the penultimate core  $\ten{g}^{(N\!-\!1)}$ is computed in a least-squares sense, as explained in Section~\ref{subsec:fixttgn}. The~pseudocode for the proposed method is outlined in Algorithm~\ref{alg:ttfw}. 
\vspace{6pt}

\begin{algorithm}[H]
	\caption{TT mode-$N$ fiber-wise}\label{alg:ttfw}
	\Indm
	\KwIn{A tensor $\ten{X} \in \R^{I_{1} \times I_{2} \times \cdots \times I_{N}}$ observed through mode-$N$ fibers,  TT rank $(R_{0}, \ldots , R_{N})$ and  row selection matrix  $\mat{s}_{r}$\\}
	\KwOut{TT cores $\ten{G}^{(1)}, \ldots, \ten{G}^{(N)}$}
	\Indp   
	\tcc{compute orthonormal bases for the ranges of the $n$th unfoldings}
	\For{$ 1 \leq n \leq N-2$}{
		Compute $\mat{a}^{(n)} \in \R^{\prod_{i = 1}^{n} I_{i} \times R_{n}}$ using the subspace constraint method, as~ discussed in Section~\ref{subsec:ssc}, or~the subspace intersection method, as~discussed in  Section~\ref{subsec:ssI}. \label{line:rangettfw}
	} 
	\tcc{compute the TT cores}
	Set $\ten{g}^{(1)} = \texttt{reshape}\left(\mat{a}^{(1)}, [1, I_1, R_1]\right)$ \;
	
	\For{$ 1 \leq n \leq N-3$}{ 
		Compute $\mat{w}^{(n)} = \mat{a}^{(n)\top} \left(\mat{a}^{(n+1)}\right)_{[1, \ldots, n ; n+1, (n+1)+1]}$\;
		Set $\ten{g}^{(n+1)} = \texttt{reshape}\left(\mat{w}^{(n)}, [R_{n}, I_{n+1}, R_{n+1}] \right)$
	}
	Compute $[\sim, \sim, \mat{v}^{(N-1)}] = \texttt{svd}\left(\mat{s}_{r}\mat{x}_{[1, \ldots,  N-1 ; N]},  R_{N-1}\right)$\;
	Set $\ten{g}^{(N)} = \texttt{reshape}\left(\mat{v}^{(N-1) \top}, [R_{N-1}, I_N, 1]\right)$\;
	Compute $\ten{g}^{(N\!-\!1)}$ by solving the linear systems in \eqref{eqn:fixttgn2}\;
\end{algorithm}

\subsubsection{Computing the Next-to-Last TT Core $\ten{g}^{(N\!-\!1)}$}
\label{subsec:fixttgn}
In the fully observed case, the~last core can be computed using the SVD of the $(N\!-\!1)$th unfolding matrix $\mat{x}_{[1,\dots,N-1;N]}$. Let its SVD be $\mat{x}_{[1,\dots,N-1;N]} = \mat{U}^{(N-1)}\,\mat{\Sigma}^{(N-1)}\,\mat{V}^{(N-1)\top}.$ Then, the~last TT core is obtained as $\ten{g}^{(N)} = \mat{\Sigma}^{(N-1)}\mat{v}^{(N-1)\top}.$  However, for~ fiber-wise observed tensors, this is not directly possible. Indeed, since $\ten{x}$ is observed through mode-$N$ fibers, some rows of $\mat{x}_{[1,\dots, N-1;N]}$ are now completely observed while other rows are entirely missing. Consequently, under~the informational-completeness condition, $\mat{v}^{(N-1)\top}$ can be estimated but $\mat{U}^{(N-1)}$ and $\mat{\Sigma}^{(N-1)}$ cannot.  To~address this issue, we set the last core $\ten{g}^{(N)} = \mat{v}^{(N-1)\top},$ i.e.,~the orthonormal basis for the observed rows $(\mat{s}_{r}\mat{x}_{[1, \ldots,  N-1 ; N]})$.  Here,  $\mat{s}_{r} \in \{0,1\}^{|\alpha| \times \prod_{i = 1}^{N-1} I_{i}}$  represents the row selection matrix for the $(N-1)$th matrix unfolding, where $|\alpha|$ denotes the number of observed mode-$N$ fibers (rows) indexed by $\alpha$. We then fix the scaling indeterminacies by computing the next-to-last TT core in a least-squares sense. For~the third-order tensor $\tilde{\ten{x}} = \ten{x}_{[1, \ldots, N-2 ; N-1; N]}$, the~slice-wise representation, using  \eqref{eqn:tt3}, can be written as
\begin{equation}
	\label{eqn:fixttgn1}
	\tilde{\mathbf{X}}_{:i:} = \mat{g}^{(<N-1)}\mat{g}^{(N-1)}_{:i:} \mat{g}^{(>N-1)}, ~ i \in \{ 1, \ldots, I_{N-1}\},
\end{equation}
where the matrices $\mat{g}^{(<N-1)} \in \R^{\prod_{i = 1}^{N-2} I_{i} \times R_{N-2}}$ and $\mat{g}^{(>N-1)} = \mat{g}^{(N)} \in \R^{R_{N-1} \times I_{N}}$ have orthonormal columns and orthonormal rows, respectively. Let $\{\mat{s}^{(i)}_{r}\}_{i=1}^{I_{N-1}}$ be the row selection matrices for the mode-2 slices of $\tilde{\ten{x}}$. The~mode-2 slices $\mat{g}^{(N\!-\!1)}_{:i:}$ of $\ten{g}^{(N-1)}$ can now be computed by solving the following linear systems:
\begin{equation}
	\label{eqn:fixttgn2}
	\resizebox{0.7\linewidth}{!}{$
		\mat{s}^{(i)}_{r}\left(\tilde{\mathbf{X}}_{:i:} \mat{g}^{(N)\top}\right)= \left(\mat{s}^{(i)}_{r}\mat{g}^{(<N-1)}\right)\mat{g}^{(N-1)}_{:i:}, ~ i \in \{ 1, \ldots, I_{N-1}\}.$
	}
\end{equation}

{
	\subsubsection{Computational Complexity}

	The computational cost of the algorithm is dominated by the subspace computation step in line~\ref{line:rangettfw} of Algorithm~\ref{alg:ttfw}. This step can be performed in parallel since the unfolding matrices are independent; hence, we report the per-processor cost. For~the $n$th unfolding matrix $\mat{x}_{[1,\ldots,n;\,n+1,\ldots,N]}$, we compute an orthonormal basis $\mat{A}^{(n)}$ for its column space via a rank-$R$ SVD of the concatenated sparse matrix $\mat{Q} \in \mathbb{R}^{I^{n} \times \left(LR + \sum_{l=1}^L (I^{n}-J_l)\right)}$, as~discussed in Section~\ref{subsec:ssI} (alternatively,  $\mat{N}$ can be used; see Section~\ref{subsec:ssc}), assuming $I_n=I$ and $R_n=R$. Here, $L$ is the number of slices used---at most $I^{N-n-1}$ (as shown in \eqref{eqn:ttfw1})---that satisfy the informational-completeness condition of Section~\ref{subsec:uniq}, and~$J_l$ is the number of observed rows in the $l$th~slice.
	
	In practice, one uses sparse (or randomized) SVD methods on $\mat{Q}$, yielding a cost of $\mathcal{O}\!\left(T\,\operatorname{nnz}(\mat{Q})\,R\right)$, where $T$ is the number of iterations and  $\operatorname{nnz}(\mat{Q})$ denotes the number of nonzeros in $\mat{Q}$ given by,  $\sum_{l=1}^L J_l R + \sum_{l=1}^L (I^{n}-J_l)$. Under~minimal oversampling, where each slice contributes as few as $J_l \sim R+c$ observed rows for some small constant $c>0$, and~all $L=I^{N-n-1}$ slices are required (worst-case), we have $\operatorname{nnz}(\mat{Q}) = L I^{n} + L(R^2 +cR-R-c) \sim  L I^{n} ~  \text{for } I^n \gg R$. Consequently, the~resulting computational cost scales as $\mathcal{O}(T I^{N\!-\!1} R)$.
}
\subsection{Uniqueness~Conditions}\label{subsec:uniq}
This section summarizes the deterministic conditions under which a  tensor $\ten{x}\in \R^{I_{1} \times I_{2} \times \cdots \times I_{N}},$  observed through mode‑$N$ fibers, admits a unique TT decomposition \(
\ten{X} = \ten{G}^{(1)} \bullet \cdots \bullet \ten{G}^{(N)}
\) (up to basis transformations (i.e., inserting $\mat{R R^{-1}}$ between consecutive TT cores, with invertible $\mat{R}$, leaves the tensor unchanged)), in~the noiseless setting. Intuitively, deterministic recovery follows from the fact that the fiber-wise observation pattern allows for an {informational overlap} to be created that algebraically locks in the solution and, consequently, reduces the completion problem to computing simple SVDs and solving a few linear~systems. \par 

To compute the TT cores $\ten{G}^{(1)},\ldots,\ten{G}^{(N\!-\!2)}$, Algorithm~\ref{alg:ttfw} utilizes orthonormal bases for the column spaces of the unfoldings
$\mat{X}_{[1,\ldots,n;\,n+1,\ldots,N]}$, $n=1,\ldots,N\!-\!2$, which exhibit a structured observation pattern, as~shown in  Figure~\ref{fig:slicestacks}. 
For each $n$, the~basis $\mat{A}^{(n)}$ is obtained via the piecewise subspace learning approach of Section~\ref{subsec:pssl}, with~the unfolding playing the role of $\mat{M}$ and the observed slice rows $\mat{s}^{(l)}_{r}\tilde{\mathbf{X}}_{:l:}$ playing the role of $\mat{m}^{(l)}_{obs}$. 
Accordingly, $\mat{A}^{(n)}$ is recovered from $\mat{N}$ (Section~\ref{subsec:ssc}) or $\mat{Q}$ (Section~\ref{subsec:ssI}). 
The last core $\ten{g}^{(N)}$ is obtained by computing the top $R_{N\!-\!1}$ right-singular vectors of the matrix formed by the observed rows of the $(N\!-\!1)$th unfolding, and~$\ten{g}^{(N\!-\!1)}$ is obtained by solving \eqref{eqn:fixttgn2}.

\begin{Theorem}[Algebraic conditions]\label{th1}
	The TT cores $\{\ten{g}^{(n)}\in\mathbb{R}^{R_{n-1}\times I_n\times R_n}\}_{n=1}^N$ are uniquely determined (up to basis transformations) if the following hold:
	\begin{itemize}
	\item[]
	\begin{itemize}
		\item[I.] \textbf{Cores \boldmath$\mathcal{G}^{(1)}\!$--$\mathcal{G}^{(N\!-\!2)}$:}
		For each $n\!=\!1,\ldots,N\!-\!2$, the~$n$th unfolding contains sufficiently many fully observed rank-$R_n$ (isorank) submatrices $\{\mat{s}^{(l)}_{r}\tilde{\mathbf{X}}_{:l:}\}$ such that 
		$\dim(\ker(\mat{N}^\top))=R_n$ (equivalently $\dim(\cap_{l=1}^L S_l)\!=\!R_n$); this is the informational-completeness condition (see \mbox{Sections~\ref{subsec:ssc} and~\ref{subsec:ssI})}.
		\item[II.] \textbf{Core \boldmath$\mathcal{G}^{(N\!-\!1)}$:}
		The systems in \eqref{eqn:fixttgn2} admit a unique solution, which is the case if
		$\mat{s}^{(i)}_{r}\mat{g}^{(<N-1)}$ has full column rank for all $i$.
		\item[III.] \textbf{Core \boldmath$\mathcal{G}^{(N)}$:} The matrix $\mat{S}_r\mat{X}_{[1,\ldots,N-1;N]}$ containing the observed mode-$N$ fibers has rank $R_{N-1}$. 
	\end{itemize}
\end{itemize}	
\end{Theorem}

\begin{Remark}
	The conditions in Theorem~\ref{th1} are not only sufficient for computation but also necessary for uniqueness. If~a core fails to satisfy its condition, it is not uniquely identifiable.
\end{Remark}
\begin{Corollary}[Generic conditions]
	Theorem~\ref{th1} holds generically~if the following hold:
	\begin{itemize}
		\item[] 
	\begin{itemize}
		\item[I'.] For each $n=1,\ldots,N\!-\!2:$
		\begin{itemize} \item[(i)] Among the observed submatrices $\{\mat{s}^{(l)}_{r} \tilde{\mathbf{X}}_{:l:}\}$ of the $n$th unfolding, some pairs of submatrices overlap in at least $R_n$ observed rows \cite{stijn2023mlsvdfsj}, Corollary~3.5. \item[(ii)] Every row of the unfolding is (partially) observed in at least one submatrix.
		\end{itemize} 
		These generically imply  Condition~I; see~\cite{stijn2023mlsvdfsj,pimentel2016dsc} for details.
		\item[II'.] Each mode-2 slice of $\ten{X}_{[1,\ldots,N-2;\,N-1;\,N]}$ contains at least $R_{N-2}$ observed rows.
		\item[III'.] At least $R_{N-1}$ mode-$N$ fibers are observed.
	\end{itemize}
\end{itemize}
\end{Corollary}

In the noiseless case, the~proposed method computes the exact TT decomposition of the fiber-wise observed tensor if the uniqueness conditions are satisfied. In~the presence of noise, each step can be computed in a least-squares sense, and~thus the~estimated  TT decomposition is expected to  be close to the true TT~decomposition.

\section{Numerical~Experiments}\label{sec:exp}
This section evaluates the performance of the proposed method across three setups. First, we experiment with synthetic data to show that our algebraic method is accurate and computationally efficient, comparing it with the TT weighted optimization (TT-WOPT)~\cite{yuan2019ttc}, {parallel matrix factorization for low‑rank TT completion (TMac‑TT)~\cite{xu2015parallel,  bengua2017silrctt},} and simple low‑rank TT completion (SiLRTC‑TT) \cite{jfeng2010svth, bengua2017silrctt}.  Second, we showcase two real-world applications, namely multidimensional harmonic retrieval (MHR) and spatiotemporal weather data completion. Finally, we show the use of the proposed algebraic approach as a ``proxy'' for efficiently performing subsequent computations. In~line with~\cite{nico2019compactrep}, we fit a non-negative CPD to the TT approximation instead of the actual~tensor. \par
\vspace{12pt}

\hspace{-24pt}
\begin{tikzpicture}
	\node[draw=black, dashed, rounded corners=5pt, align=left, inner sep=5pt, fill=gray!5]
	at (0,0) {%
		\begin{minipage}{0.96\textwidth}
			{\small 
				Baselines:
				\begin{itemize}
					\item TT-WOPT is an error minimization method that fits the parameters of a fixed rank TT model by minimizing a weighted least-squares loss using the gradient descent approach; see~\cite{yuan2019ttc} for~details.  
					
					\item{TMac-TT is an error minimization method that enforces low TT rank by fitting matrix factorizations to each mode unfolding of the tensor and updating the factors using alternating least squares; see~\cite{xu2015parallel, bengua2017silrctt} for more details.}
					
					\item  SiLRTC-TT is a rank minimization method  that minimizes a nuclear norm relaxation of the TT rank of the tensor. The~resulting convex problem is solved by a singular value thresholding algorithm; see~\cite{jfeng2010svth, bengua2017silrctt} for more~details.  
					
				\end{itemize}
			}
		\end{minipage}
	};
\end{tikzpicture}
\vspace{6pt}

The experiments were performed on an HP EliteBook 845 G8 Notebook PC with an AMD Ryzen 7 PRO 5850U CPU and 32 GB~RAM.

\subsection{Synthetic~Data}
\subsubsection{Completion from Noisy~Observations}
\begin{figure}[H]
	\centering
	\begin{tikzpicture}
		\begin{groupplot}[
			group style={group size=2 by 1, horizontal sep=2.4cm},
			width=0.4\textwidth,
			legend to name=sharedlegendepx1,
			legend style={draw=none},
			legend columns=4,
			]
			\nextgroupplot[
			xlabel={SNR},
			ylabel={Relative Error},
			axis lines=middle,
			enlargelimits=true,
			axis x line = bottom,
			axis y line =left,
			grid=both,
			grid style={line width=0.2pt, draw=gray!70},
			xtick={-10, 0,10,20,30,40,50},
			ymode=log,
			ymin=0,
			ymax=4,
			scale=1,
			every axis x label/.append style={
				at={(axis description cs:0.5,-0.3)},
				anchor=north, color=gray
			},
			every axis y label/.append style={
				at={(axis description cs:-0.32,0.5)},
				anchor=south,
				rotate=90, color=gray},
			x axis line style={yshift=-10pt,-,color=gray},
			xticklabel style={yshift=-10pt,color=gray},
			y axis line style={xshift=-10pt,-,color=gray},
			yticklabel style={xshift=-10pt,color=gray},
			xtick style={yshift=-10pt, line width=0.5pt},
			ytick style={xshift=-10pt, line width=0.5pt},
			xtick align=inside,
			ytick align=inside,
			]
			
			\addplot[thick, color={rgb,255:red,139; green,0; blue,139}, mark=square*, mark size=2pt] coordinates {(-10, 0.9234)(0,0.2900) (10,0.0868) (20, 0.0391) (30, 0.0187) (40, 0.0088) (50, 0.0047)};
			\addlegendentry{TT mode-$N$ fiber-wise}
			
			\addplot[thick, blue, mark=*, mark size=2pt] coordinates {(-10,  0.3832)(0, 0.1130) (10, 0.0356) (20, 0.0113) (30, 0.0036) (40,0.0011) (50,0.0005)};
			\addlegendentry{TT-WOPT}
			
			\addplot[thick, color=black, mark=triangle*, mark size=2pt] coordinates {(-10, 3.0838) (0, 0.9295) (10, 0.2525) (20, 0.0468) (30, 0.0167) (40, 0.0138) (50,0.0118)};
			\addlegendentry{SiLRTC-TT}
			
			\addplot[thick, color=olive, mark=diamond*, mark size=2pt] coordinates {(-10, 3.1516) (0, 0.7397) (10, 0.1537) (20, 0.0405) (30, 0.0114) (40, 0.0035) (50,0.0019)};
			\addlegendentry{TMac-TT}
			
			\nextgroupplot[
			xlabel={SNR},
			xtick={-10, 0,10,20,30,40,50},
			ylabel={\textcolor{gray}{Time (s)}},
			axis lines=middle,
			enlargelimits=true,
			axis x line = bottom,
			axis y line = left,
			ymode=log,
			ymin=0.01,
			ymax=10,
			grid=both,
			grid style={line width=0.2pt, draw=gray!70},
			every axis x label/.append style={
				at={(axis description cs:0.5,-0.3)},
				anchor=north, color=gray
			},
			every axis y label/.append style={
				at={(axis description cs:-0.27,0.5)},
				anchor=south,
				rotate=90, color=gray},
			x axis line style={yshift=-10pt,-,color=gray},
			xticklabel style={yshift=-10pt,color=gray},
			y axis line style={xshift=-10pt,-,color=gray},
			yticklabel style={xshift=-10pt,color=gray},
			xtick style={yshift=-10pt, line width=0.5pt},
			ytick style={xshift=-10pt, line width=0.5pt},
			xtick align=inside,
			ytick align=inside,
			]

			\addplot[thick, color={rgb,255:red,139; green,0; blue,139}, mark=square*, mark size=2pt] coordinates {(-10, 0.0293)(0,0.0243) (10, 0.0237) (20, 0.0223) (30, 0.0230) (40,0.0223) (50, 0.0222)};
			\addlegendentry{TT mode-$N$ fiber-wise}
			
			\addplot[thick, blue, mark=*, mark size=2pt] coordinates {(-10, 1.2053)(0, 1.1922) (10, 1.1565) (20, 1.1848) (30, 1.1486) (40, 1.1436) (50, 1.1145)};
			\addlegendentry{TT-WOPT}
			
			\addplot[thick, color=black, mark=triangle*, mark size=2pt] coordinates {(-10, 8.1436) (0, 8.0162) (10, 4.6384) (20, 2.5487) (30, 1.7886) (40, 1.7599) (50, 1.8133)};
			\addlegendentry{SiLRTC-TT}
			
			\addplot[thick, color=olive, mark=diamond*, mark size=2pt] coordinates {(-10, 0.9421) (0, 0.9785) (10,0.8510) (20, 0.8077) (30, 0.7835) (40, 0.7917) (50, 0.7970)};
			\addlegendentry{TMac-TT}
			
		\end{groupplot}

		\node at ($ (group c1r1.south)!0.5!(group c2r1.south) + (0,-1.9cm) $) 
		{\pgfplotslegendfromname{sharedlegendepx1}};
	\end{tikzpicture}
	\caption{(\textbf{Left}) The accuracy of our method is slightly lower than that of TT-WOPT and TMac-TT, which is expected due to its reliance on only standard NLA operations. (\textbf{Right}) However, it provides a significant computational~speedup.}
	\label{fig:synexp}
\end{figure}
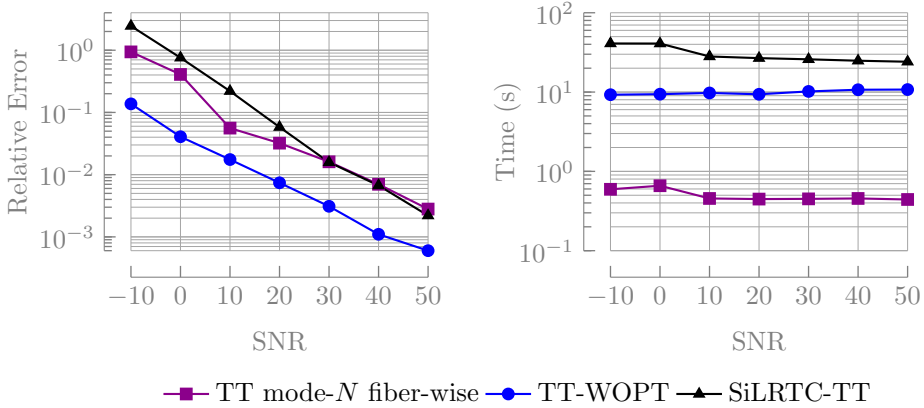

A 4th-order tensor in TT format $\ten{x} \in \R^{15  \times 15 \times 15 \times 15}$ is generated with $\text{rank}_{TT}(\ten{x}) = (1, 3, 3, 4, 1)$ by sampling entries of the TT cores from a standard normal distribution. Random Gaussian i.i.d. noise \(\ten{n}\) is added to the generated tensor, resulting in a noisy tensor \(\ten{x}_{noisy} = \ten{x} + \ten{n}\) with a fixed signal-to-noise ratio (SNR), defined as \(\text{SNR} = 20 \log_{10} \frac{\| \ten{x}\|}{\| \ten{n}\|}\). Forty percent of the mode-4 fibers are randomly removed from the noisy tensor. Algorithm~\ref{alg:ttfw} is run along with the reference algorithms to compute the low-rank TT approximations with SNR varying from $-10$ dB to $50$ dB. The~reference algorithms are run with default parameters (see~\cite{yuan2019ttc, bengua2017silrctt}), except~for parameter $\alpha$ in SiLRTC-TT, which is determined through an extensive grid search. The~median of the relative error, defined as \(\text{Relative Error} = \frac{\|\ten{x} - \hat{\ten{x}}\|_F}{\|\ten{x}\|_F}\), is plotted across 30 trials in Figure~\ref{fig:synexp} (left). The~relative error of the proposed algorithm is observed to be lower than that of the SiLRTC‑TT algorithm, while the lowest error is noted for TT‑WOPT. {TMac‑TT shows marginally lower accuracy than our approach in high‑noise regimes, but~marginally higher accuracy in low‑noise regimes. TT‑WOPT and TMac‑TT perform (marginally) better because they explicitly aim to minimize the error.}  In comparison, the~accuracy of the proposed algorithm is very good despite the fact that it only relies on standard NLA operations. The~time required to compute the approximation is shown in Figure~\ref{fig:synexp} (right), which shows that the proposed method is more than a magnitude faster, and~the effect of the SNR on computation time is not~significant. 
\unskip
\subsubsection{Scalability}
\begin{figure}[H]
	\centering
	\begin{tikzpicture}
		\begin{groupplot}[
			group style={group size=2 by 1, horizontal sep=1.8cm},
			width=0.42\textwidth,
			legend to name=sharedlegendexp2,
			legend style={draw=none},
			legend columns=4,
			]
			\nextgroupplot[
			xlabel={I},
			ylabel={Relative Error},
			axis lines=middle,
			enlargelimits=true,
			axis x line = bottom,
			axis y line =left,
			grid=both,
			grid style={line width=0.2pt, draw=gray!70},
			xtick={10,20,30,40,50},
			ymode=log,
			ymin=0,
			ymax=0.5,
			scale=0.9,
			every axis x label/.append style={
				at={(axis description cs:0.5,-0.3)},
				anchor=north, color=gray
			},
			every axis y label/.append style={
				at={(axis description cs:-0.32,0.5)},
				anchor=south,
				rotate=90, color=gray},
			x axis line style={yshift=-10pt,-,color=gray},
			xticklabel style={yshift=-10pt,color=gray},
			y axis line style={xshift=-10pt,-,color=gray},
			yticklabel style={xshift=-10pt,color=gray},
			xtick style={yshift=-10pt, line width=0.5pt},
			ytick style={xshift=-10pt, line width=0.5pt},
			xtick align=inside,
			ytick align=inside,
			]
			
			\addplot[thick, color={rgb,255:red,139; green,0; blue,139}, mark=square*, mark size=2pt] coordinates {(10,0.0556292146361901) (20,0.0338537704796569) (30,0.0184026422733671) (40,0.0118199061005346) (50, 0.00864811225339755)};
			
			\addplot[thick, blue, mark=*, mark size=2pt, mark options={solid, fill=blue}] coordinates {(10,0.0134817459882531) (20,0.00478483707661231) (30, 0.00261745219910095) (40,0.00171246829151422) (50, 0.00123645164138533)};

			\addplot[thick, color=black, mark=triangle*, mark size=2pt] coordinates {(10, 0.0678043949129817) (20,0.038704601682288) (30,0.0235621296631998) (40,0.0201493398208418) (50,0.0171493398208418)};
			
			\addplot[thick, color=olive, mark=diamond*, mark size=2pt] coordinates {(10, 0.0315290567501486) (20,0.0167785636698569) (30,0.00925840173454683) (40,0.00762723475827298) (50,0.00614991443334326)};
			
			\nextgroupplot[
			xlabel={I},
			xtick={10,20,30,40,50},
			ylabel={\textcolor{gray}{Time (s)}},
			ytick={0, 0.1, 10e0, 10e1, 10e2, 10e3},
			axis lines=middle,
			enlargelimits=true,
			axis x line = bottom,
			axis y line = left,
			ymode=log,
			ymin=0,
			ymax=10000,
			xmax=50.5,
			scale=0.9,
			grid=both,
			grid style={line width=0.2pt, draw=gray!70},
			every axis x label/.append style={
				at={(axis description cs:0.5,-0.25)},
				anchor=north, color=gray
			},
			every axis y label/.append style={
				at={(axis description cs:-0.23,0.5)},
				anchor=south,
				rotate=90, color=gray},
			x axis line style={yshift=-10pt,-,color=gray},
			xticklabel style={yshift=-10pt,color=gray},
			y axis line style={xshift=-10pt,-,color=gray},
			yticklabel style={xshift=-10pt,color=gray},
			xtick style={yshift=-10pt, line width=0.5pt},
			ytick style={xshift=-10pt, line width=0.5pt},
			xtick align=inside,
			ytick align=inside,
			minor y tick num=9,
			]

			\addplot[thick, color={rgb,255:red,139; green,0; blue,139}, mark=square*, mark size=2pt] coordinates {(10, 0.0099991) (20, 0.0558603) (30, 0.2633652 ) (40, 0.7074924) (50, 1.228323)};
			\addlegendentry{TT mode-$N$ fiber-wise}
			
			\addplot[thick, blue, mark=*, mark size=2pt, mark options={solid, fill=blue}] coordinates {(10,1.354057) (20, 2.5987039) (30, 5.9689655) (40, 27.1947772) (50, 88.8138599)};
			\addlegendentry{TT-WOPT}
			
			\addplot[thick, color=black, mark=triangle*, mark size=2pt] coordinates {(10,1.2129438) (20, 7.4644277) (30,43.8713847) (40, 469.3565242) (50, 2799.4760918)};
			\addlegendentry{SiLRTC-TT}
			
			\addplot[thick, color=olive, mark=diamond*, mark size=2pt] coordinates {(10,0.2367545) (20, 3.2415373) (30,9.5550967) (40, 57.585497) (50, 156.9289721)};
			\addlegendentry{TMac-TT}
			
			
			\coordinate (Aa) at (axis cs:30, 43.8713847);
			\coordinate (Ab) at (axis cs:30, 2799.4760918);
			\coordinate (Ac) at (axis cs:50, 2799.4760918);
			\draw[thin, dashed, color=black] (Aa) -- (Ab); 
			\draw[thin, dashed, color=black] (Ab) -- (Ac) node[midway, xshift=0.5mm, yshift=-3mm,  left] {\small $\times 57.3$};
			
			\coordinate (Ba) at (axis cs:30, 5.9689655);
			\coordinate (Bb) at (axis cs:50, 5.9689655);
			\coordinate (Bc) at (axis cs:50, 88.8138599);
			\draw[thin, dashed, color=blue] (Ba) -- (Bb);
			\draw[thin, dashed, color=blue] (Bb) -- (Bc) node[midway, xshift=0, yshift=-1.5mm,  left] {\small $\times 14.9$};
			
			\coordinate (Ca) at (axis cs:30, 0.2633652);
			\coordinate (Cb) at (axis cs:50, 0.2633652);
			\coordinate (Cc) at (axis cs:50, 1.228323);
			\draw[ thin, dashed, color={rgb,255:red,139; green,0; blue,139}] (Ca) -- (Cb);
			\draw[thin, dashed, color={rgb,255:red,139; green,0; blue,139}] (Cb) -- (Cc) node[midway, xshift=0, yshift=-0.5mm,  left] { \small $\times 4.7$};

			\coordinate (Aa) at (axis cs:30, 9.5550967);
			\coordinate (Ab) at (axis cs:30, 156.9289721);
			\coordinate (Ac) at (axis cs:50, 156.9289721);
			\draw[thin, dashed, color=olive] (Aa) -- (Ab); 
			\draw[thin, dashed, color=olive] (Ab) -- (Ac) node[midway, xshift=10mm, yshift=2mm,  left] {\small $\times 16.4$};

		\end{groupplot}

		\node at ($ (group c1r1.south)!0.5!(group c2r1.south) + (0,-1.9cm) $) 
		{\pgfplotslegendfromname{sharedlegendexp2}};

	\end{tikzpicture}
	\caption{(\textbf{Left}) Accuracy improves as $I$ increases since relatively more data becomes available per parameter (with the missing rate and TT rank held constant). TT-WOPT achieves the highest accuracy. (\textbf{Right}) Meanwhile, our method outperforms all reference algorithms in terms of computation~time.}
	\label{fig:synscalability}
\end{figure}
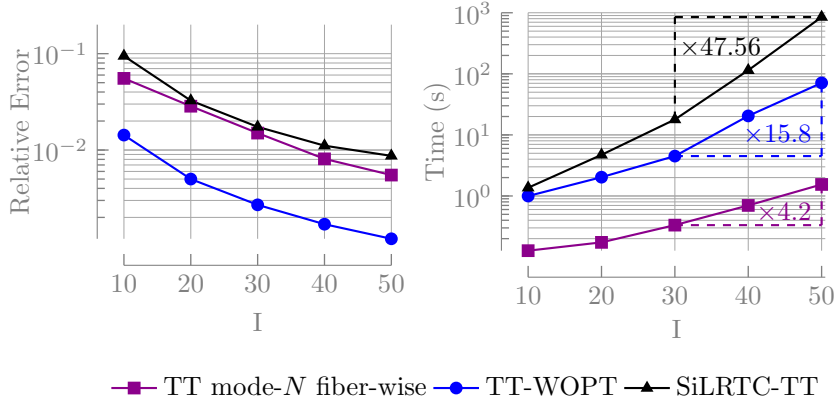

This experiment compares the scalability of the proposed approach with the reference algorithms. The~reference algorithms are based on optimization and are terminated when the relative change in the function value falls below $10^{-7}$. If~this level of precision is not achieved, the~algorithms  terminate upon reaching the maximum number of iterations (maximum iterations  \;=\; $7 \times 10^2$).  A~$4$th-order tensor in TT format $\ten{x} \in \R^{I \times I \times I \times I}$ is generated with a fixed $\text{rank}_{TT}(\ten{x}) = (1,  4, 4, 4, 1)$, where $I$ is varied from $10$ to $50$. Random Gaussian i.i.d. noise is added to the generated tensor to achieve an SNR of $25$ dB, and~$35$\% of the mode-$4$ fibers are randomly removed.  Figure~\ref{fig:synscalability} (left) shows the median relative error computed over 25 trials and plotted as a function of $I$. The~TT-WOPT algorithm is shown to achieve the highest accuracy, while a good but not optimal accuracy is achieved by our method. Notably, for~the proposed approach, accuracy consistently improves as the problem size increases. This behaviour is expected because, for~a fixed missing rate and TT rank, an~increase in \(I\) results in relatively more observed data being used to estimate the model parameters (TT cores). A similar reasoning applies to the reference optimization methods; see~\cite{jfeng2010svth, bengua2017silrctt} for more details.   In~Figure~\ref{fig:synscalability} (right), it is shown that the computation time for the proposed method is the lowest, in~line with expectations. Specifically, in~the current experimental settings, as~$I$ increases from $30$ to $50$, {the time increases by a factor of $4.7$ in our approach, while the time increases by factors of $14.9,$ $16.4$ and $57.3$ for TT-WOPT, TMac-TT} and SiLRTC-TT algorithms, respectively. 

\subsection{Real-Life~Applications}
\unskip
\subsubsection{Multidimensional Harmonic~Retrieval}
Harmonic retrieval is a classical problem in signal processing. The~goal is to estimate the parameters of a signal that is modelled as a sum of complex exponentials. MHR is the natural multidimensional extension. An~MHR data tensor sampled over the $(D+1)$th-order tensor grid can be modelled as
\begin{equation}
	\label{eqn:mhr}
	x_{i_1 i_2 \cdots i_d k}=\sum_{r=1}^R s_r(k) \prod_{d=1}^D \mathrm{e}^{j\left(i_d-1\right) \mu_r^{(d)}}+n_{i_1 i_2 \cdots i_D k},
\end{equation}
where $j^2=-1$ and $s_r(k)$ is the $k$th complex symbol carried by the $r$th multidimensional harmonic. The~noise $n_{i_1 i_2 \cdots i_D k}$ is modelled as zero-mean i.i.d additive Gaussian~noise. \par

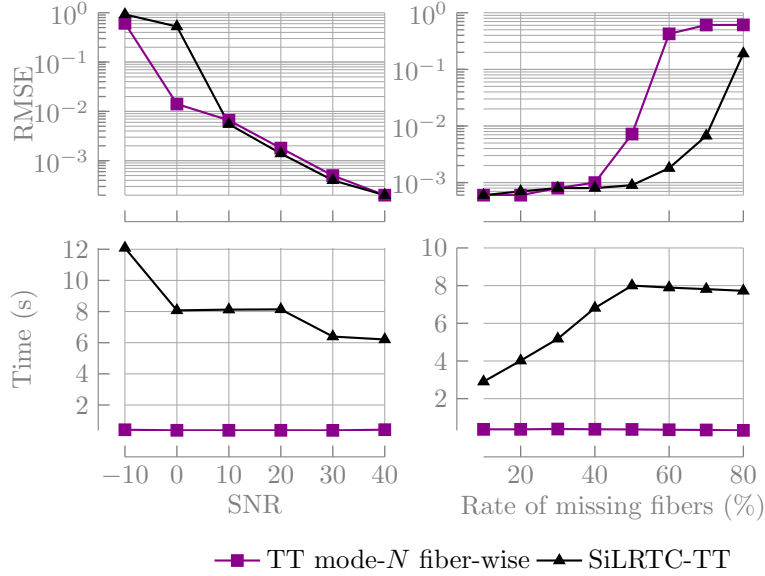
\begin{figure}[H]
	\centering
	\begin{tikzpicture}
		\begin{groupplot}[
			group style={
				group size=2 by 2,
				horizontal sep=1.3cm,
				vertical sep=0.7cm,
				group name=mdhrgroup
			},
			width=5cm,
			height=4cm,
			legend to name=sharedlegendexp3,
			legend style={draw=none},
			legend columns=-1,
			]
			
			\nextgroupplot[
			xlabel={},
			ylabel={RMSE },
			axis lines=middle,
			enlargelimits=true,
			axis x line = bottom,
			axis y line =left,
			grid=both,
			grid style={line width=0.2pt, draw=gray!70},
			ymode=log,
			ymin=0,
			ymax=1,
			xtick={-10, 0,10,20,30,40},
			xticklabels={},
			every axis x label/.append style={
				at={(axis description cs:0.5,-0.3)},
				anchor=north, color=gray
			},
			every axis y label/.append style={
				at={(axis description cs:-0.32,0.5)},
				anchor=south,
				rotate=90, color=gray},
			x axis line style={yshift=-10pt,-,color=gray},
			xticklabel style={yshift=-10pt,color=gray},
			y axis line style={xshift=-10pt,-,color=gray},
			yticklabel style={xshift=-10pt,color=gray},
			xtick style={yshift=-10pt, line width=0.5pt},
			ytick style={xshift=-10pt, line width=0.5pt, color=gray},
			xtick align=inside,
			ytick align=inside,
			]
			\addplot[thick, color={rgb,255:red,139; green,0; blue,139}, mark=square*, mark size=2pt] coordinates {(-10, 0.6099)(0,0.0141) (10, 0.0066) (20,0.0018) (30,0.0005) (40,0.0002)};
			
			\addplot[thick, color=black, mark=triangle*, mark size=2pt] coordinates {(-10, 0.9249) (0, 0.5261) (10, 0.0055) (20, 0.0014) (30, 0.0004) (40, 0.0002)};
			
			\nextgroupplot[
			xlabel={},
			ylabel={},
			axis lines=middle,
			enlargelimits=true,
			axis x line = bottom,
			axis y line =left,
			grid=both,
			grid style={line width=0.2pt, draw=gray!70},
			ymode=log,
			ymin=0,
			ymax=1,
			xtick={20,40,60,80},
			xticklabels={},
			every axis x label/.append style={
				at={(axis description cs:0.5,-0.3)},
				anchor=north, color=gray
			},
			every axis y label/.append style={
				at={(axis description cs:-0.3,0.5)},
				anchor=south,
				rotate=90, color=gray},
			x axis line style={yshift=-10pt,-,color=gray},
			xticklabel style={yshift=-10pt,color=gray},
			y axis line style={xshift=-10pt,-,color=gray},
			yticklabel style={xshift=-10pt,color=gray},
			xtick style={yshift=-10pt, line width=0.5pt},
			ytick style={xshift=-10pt, line width=0.5pt, color=gray},
			xtick align=inside,
			ytick align=inside,
			]
			
			\addplot[thick, color={rgb,255:red,139; green,0; blue,139}, mark=square*, mark size=2pt] coordinates {(10, 0.0006) (20, 0.0006) (30,0.0008) (40, 0.0010) (50, 0.0072) (60, 0.4241)(70, 0.6080) (80, 0.6100)};
			
			\addplot[thick, color=black, mark=triangle*, mark size=2pt] coordinates {(10,0.0006) (20, 0.0007) (30,0.0008) (40, 0.0008) (50, 0.0009) (60, 0.0018)
				(70, 0.0067) (80, 0.190)};
			
			\nextgroupplot[
			xlabel={SNR},
			ylabel={Time (s) },
			axis lines=middle,
			enlargelimits=true,
			axis x line = bottom,
			axis y line =left,
			grid=both,
			grid style={line width=0.2pt, draw=gray!70},
			ytick={2, 4, 6, 8, 10, 12},
			xtick={-10, 0, 10,20,30,40},
			every axis x label/.append style={
				at={(axis description cs:0.5,-0.3)},
				anchor=north, color=gray
			},
			every axis y label/.append style={
				at={(axis description cs:-0.3,0.5)},
				anchor=south,
				rotate=90, color=gray},
			x axis line style={yshift=-10pt,-,color=gray},
			xticklabel style={yshift=-10pt,color=gray},
			y axis line style={xshift=-10pt,-,color=gray},
			yticklabel style={xshift=-10pt,color=gray},
			xtick style={yshift=-10pt, line width=0.5pt},
			ytick style={xshift=-10pt, line width=0.5pt, color=gray},
			xtick align=inside,
			ytick align=inside,
			]
			
			\addplot[thick, color={rgb,255:red,139; green,0; blue,139},  mark=square*, mark size=2pt, mark options={solid, fill={rgb,255:red,139; green,0; blue,139}}] coordinates {(-10, 0.4110) (0,0.3754) (10,0.3753) (20,0.3763) (30,0.3732) (40, 0.4155)};
			
			\addplot[thick, color=black, mark=triangle*, mark size=2pt, mark options={solid, fill=black}] coordinates {(-10, 12.0838) (0, 8.0754) (10, 8.1237) (20, 8.1403) (30, 6.3955) (40, 6.2069)};
			
			\nextgroupplot[
			xlabel={Rate of missing fibers (\%)},
			ylabel={},
			ymax=10,
			ytick={2, 4, 6, 8, 10},
			xtick={20,40,60,80},
			axis lines=middle,
			enlargelimits=true,
			axis x line = bottom,
			axis y line =left,
			grid=both,
			grid style={line width=0.2pt, draw=gray!70},
			every axis x label/.append style={
				at={(axis description cs:0.5,-0.3)},
				anchor=north, color=gray
			},
			every axis y label/.append style={
				at={(axis description cs:-0.3,0.5)},
				anchor=south,
				rotate=90, color=gray},
			x axis line style={yshift=-10pt,-,color=gray},
			xticklabel style={yshift=-10pt,color=gray},
			y axis line style={xshift=-10pt,-,color=gray},
			yticklabel style={xshift=-10pt,color=gray},
			xtick style={yshift=-10pt, line width=0.5pt},
			ytick style={xshift=-10pt, line width=0.5pt, color=gray},
			xtick align=inside,
			ytick align=inside,]
			
			\addplot[thick,  color={rgb,255:red,139; green,0; blue,139},  mark=square*, mark size=2pt, mark options={solid, fill={rgb,255:red,139; green,0; blue,139}}] coordinates {(10, 0.3711) (20, 0.3733) (30, 0.3882) (40, 0.3759) (50, 0.3668) (60,0.3501)(70,0.3390) (80, 0.3183)};
			\addlegendentry{TT mode-$N$ fiber-wise}
			
			\addplot[thick, color=black, mark=triangle*, mark size=2pt, mark options={solid, fill=black}] coordinates {(10, 2.9044) (20, 4.0179) (30, 5.1800) (40, 6.8133) (50, 8.0035) (60, 7.8951)	(70, 7.8145) (80, 7.7241)};
			\addlegendentry{SiLRTC-TT}
		\end{groupplot}
		
		\node[anchor=north, font=\normalsize] at ($(mdhrgroup c1r2.outer south)!0.5!(mdhrgroup c2r2.outer south) + (1cm, -0.05cm)$) {\pgfplotslegendfromname{sharedlegendexp3}};
	\end{tikzpicture}
	\caption{(\textbf{Top left}) The proposed method achieves lower RMSE at low SNRs and performs comparably to SiLRTC-TT at higher SNRs. (\textbf{Bottom left}) It also requires less computation time, while SiLRTC-TT converges slowly at low SNRs. (\textbf{Top right}) Our approach remains accurate even with up to a 50\% missing rate, but~beyond that, performance degrades sharply as its working conditions no longer hold. (\textbf{Bottom right}) The proposed method is faster overall, while SiLRTC-TT slows down at high missing rates due to increased optimization~challenges. }
	\label{fig:mdhexp}
\end{figure}

A $5$th-order data tensor $\ten{x} \in \C^{10 \times 10 \times 10 \times 10 \times 25}$ is generated by a CPD model of rank $R = 4$, using \eqref{eqn:mhr}. The~parameter $D$ is set to $4$, and~binary phase shift keying sources $s_r(k) \in \{-1, 1\}$ of  length $K = 25$ are used. The~parameter vectors are set as follows: $\boldsymbol{\mu}^{(1)}=[1, -0.5, 0.1, -0.8],$ $\boldsymbol{\mu}^{(2)}=[-0.5, 1, -0.9, 0.2],$ $\boldsymbol{\mu}^{(3)}=[-0.2, -0.6, 1.0, 0.4],$ and $\boldsymbol{\mu}^{(4)}=[-0.8, 0.4, 0.3, -0.1]$.  The~parameter settings and data generation are performed as described in~\cite{nico2014curseofdim}.  In~the first experiment,  Gaussian noise is added to the generated tensor with SNR varying from -10 dB to 40 dB. Forty percent of fibers in mode-$5$ are randomly removed. The~estimated TT approximations are subsequently used to compute the parameter vectors $\hat{\boldsymbol{\mu}}^{(d)}$ using the classical ESPRIT algorithm~\cite{roy1989esprit}. Root mean square error (RMSE), which is defined as $ \text{RMSE}=\sqrt{\frac{1}{R D} \sum_{r=1}^R \sum_{d=1}^D\left(\mu_r^{(d)}-\hat{\mu}_r^{(d)}\right)^2},$ is used to assess the accuracy. The~results are summarized in the left column of Figure~\ref{fig:mdhexp}. Next,  the~SNR is fixed to 25 dB, and~the effect of the missing rate on the accuracy is studied. The~RMSE and computation time are plotted as a function of the missing fiber rate in the right column of Figure~\ref{fig:mdhexp}.

\subsubsection{Spatiotemporal Weather Data~Imputation}\label{stf}

Spatiotemporal weather data are typically collected at fixed spatial coordinates, with~measurements at each location varying over time. In~practice, however, it is often not feasible to gather (or store) data at every location---especially when the spatial resolution is high. In~such cases, time series data are recorded only for a subset of spatial coordinates across selected time windows. These data may be organized in a tensor of size---for example, ``longitudes $\times$ latitudes  $\times$  year  $\times$ day of year''---with observations in a fiber-wise pattern. In~this experiment, such a dataset, consisting of the maximum temperature time series (TMAX in $^\circ$C) from the NASA POWER database (these data were obtained from the NASA Langley Research Center (LaRC) POWER Project funded through the NASA Earth Science/Applied Science Program: \url{https://power.larc.nasa.gov})  
is used. The~dataset comprises \(5\,478\) daily observations from 01~January~2005 through 31~December~2019, and~spans the region bounded by \(4.0^\circ\mathrm{E}\) to \(50.5^\circ\mathrm{E}\) longitude and  \(30.0^\circ\mathrm{N}\) to \(54.5^\circ\mathrm{N}\) latitude, on~a regular \(0.5^\circ \times 0.5^\circ\) grid. The~code to download the dataset is provided in~\cite{sofi2024stf}. The~data are reshaped into a $4$th-order tensor, organized yearly, with~each year comprising 366 days  (day-of-year alignment is performed, with~missing entries imputed via nearest neighbor averaging), resulting in a shape of $94 \times 50 \times 15 \times 366$. We simulate a scenario in which time series data are observed only for a subset of spatial coordinates during specific years. This results in a spatiotemporal tensor whose fibers along the last mode are either fully observed or entirely~missing.  \par

In the first experiment, the~mode-$4$ fibers of the data tensor are randomly removed, with~the rate of missing fibers varying from $40\%$ to $65\%$. The~approximation is computed using our method for different TT ranks, and~the median relative error between the ground truth and the estimated tensor (i.e., the~overall error, which includes both prediction error and reconstruction  error) is recorded over 30 trials.   In~Table~\ref{tab:expstf}, it is observed that the approximation is improved by increasing the TT rank. A~reasonably good approximation is obtained even when up to 65\% of the mode-4 fibers are completely missing. However, the~error is observed to rise once the missing fibers exceed this rate. A~sharp rise in relative error is noted when the approximation rank is increased under a high missing rate (see, e.g.,~ the  value highlighted in Table~\ref{tab:expstf}). With~such a high TT rank $\left (\text{i.e., } \operatorname{rank}_{TT}(\hat{\ten{x}}) = (1, 42, 42, 46, 1) \right)$, the~overlapping conditions are no longer satisfied, resulting in the estimated approximation being rendered inaccurate; see Section~\ref{subsec:uniq}. Nevertheless, a~valid approximation can still be obtained if the TT rank is low enough to satisfy the recovery~conditions.

\setlength{\tabcolsep}{12pt}
\begin{table}[ht]
	\centering
	\sisetup{table-format=1.3}
	\caption{For a fixed missing rate, the~completion accuracy improves with increasing TT rank. For~a fixed TT rank, the~error increases only slightly, showing that the method can achieve reasonably good approximations even at higher missing rates if the working conditions are met. If~not (as in the highlighted case), the~approximation becomes~inaccurate.}
	\label{tab:expstf}
	\begin{tabularx}{\textwidth}{l l *{5}{S[table-format=1.3]}}
		\toprule
		\multicolumn{2}{c}{} & \multicolumn{5}{c}{\text{Rate of missing fibers (\%)}} \\
		\cmidrule(lr){2-7}
		\text{$\operatorname{rank}_{TT}(\hat{\ten{x}})$} &  {40} & {45} & {50} & {55} & {60} & {65} \\
		\midrule
		\text{(1, 10, 10, 14, 1)}       
		& 0.1393 & 0.1395 & 0.1396 & 0.1396 & 0.1396 & 0.1397 \\
		\addlinespace
		\text{(1, 18, 18, 22, 1)}  
		& 0.1174 & 0.1176 & 0.1177 & 0.1180 & 0.1181 & 0.1180 \\
		\addlinespace
		\text{(1, 26, 26, 30, 1)}  
		& 0.1044 & 0.1044 & 0.1045 & 0.1046 & 0.1049 & 0.1050 \\
		\addlinespace
		\text{(1, 34, 34, 38, 1)} 
		& 0.0943 & 0.0944 & 0.0946 & 0.0947 & 0.0956 & 0.0956 \\
		\addlinespace
		\text{(1, 42, 42, 46, 1)} 
		& 0.0869 & 0.0869 & 0.0874 & 0.0873 & 0.0881 & \textbf{0.6304} \\
		\bottomrule
	\end{tabularx}
\end{table}

\vspace{-6pt}

\begin{figure}[H]
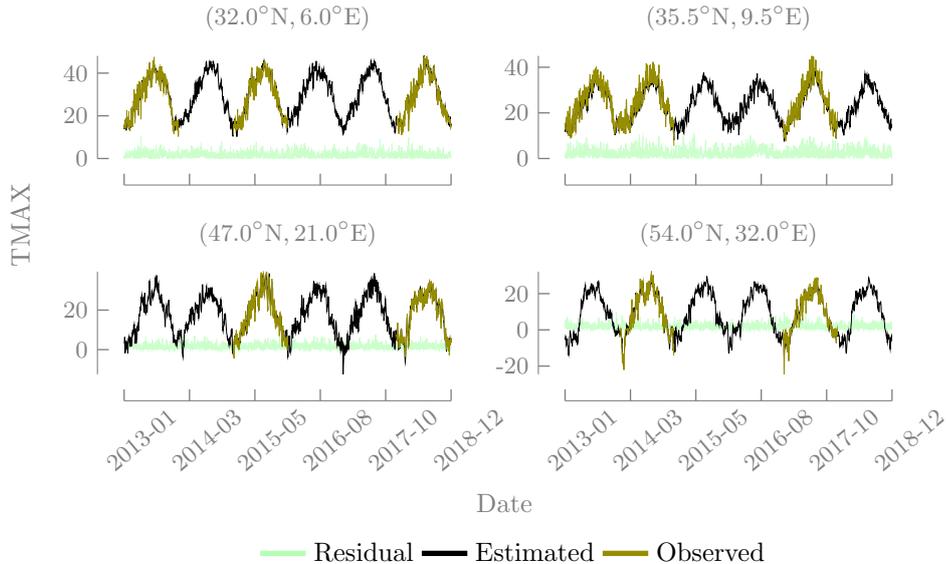

	\centering

	\caption{The estimated time series closely matches the observations. Even when data are completely missing for consecutive years, our estimation remains reasonably accurate, as~evidenced by the small residuals. This is due to the fact that the data exhibit a low-rank structure and easily identifiable~patterns.}
	\label{fig:stcts}
\end{figure}

Next, the~rate of missing fibers is set to 50\%, and~mode-4 fibers are randomly removed. The~approximation is computed using \(\operatorname{rank}_{TT}(\hat{\ten{x}}) = (1, 15, 15, 50, 1)\), which is determined by analyzing the singular values of the matrix unfoldings. The~median relative error across 10 trials is found to be below 9.7\%. Figure~\ref{fig:stcts} visualizes the observed, estimated, and~residual values (i.e., the~absolute differences between the ground truth and the estimates) for segments of the time series between 01~January~2013 and 12~December~2018 at four locations. It is observed that the maximum temperature (TMAX) exhibits a sinusoidal pattern, with~peaks corresponding to summer temperatures and valleys corresponding to winter~temperatures.

\unskip

\subsection{TT Approximation as a Prior for Efficient~Computation}
\label{proxy}

One can use our algebraic method as an initialization for optimization-based methods, improving computational speed and potentially reducing the risk of convergence to local minima, since the algebraic method---which works under deterministic conditions---yields a solution already close to the true one. Moreover, the~solution obtained from the proposed method can also be directly used for downstream tasks, particularly in low-noise settings. In~line with~\cite{nico2019compactrep}, our experiments show that the estimated TT approximation can also serve as a proxy for efficiently computing other tensor~decompositions.

\subsubsection{Initialization of Optimization~Methods}
\label{proxyinit}

In a noiseless setting, the~proposed algorithm computes the exact TT decomposition when the uniqueness conditions are satisfied. However, in~noisy settings, the~approximation is good but not optimal; further refinement can be performed using optimization methods. This experiment compares TT-WOPT---run with different \emph{random} initializations---to a hybrid approach in which the \emph{algebraic} method serves as an initialization for~TT-WOPT.    

\begin{figure}[H]
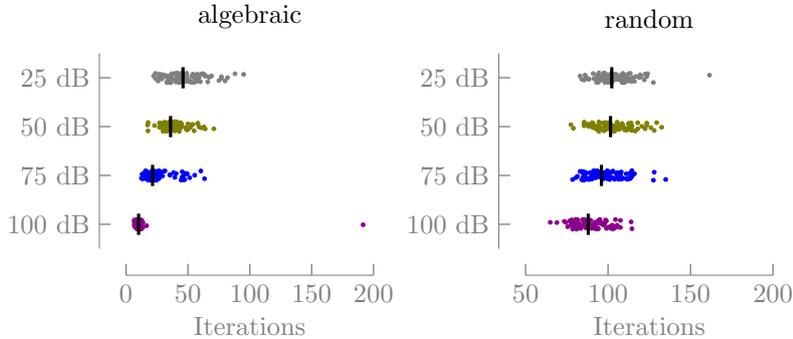

	\definecolor{mycolor1}{rgb}{0.545,0,0.545}%
	\definecolor{mycolor2}{rgb}{0,0,1}%
	\definecolor{mycolor3}{rgb}{0.502, 0.502, 0}%
	\definecolor{mycolor4}{rgb}{0.5,0.5,0.5}%
	
	\centering

	\caption{Both strategies achieve comparable accuracy, but~with algebraic initialization, the~TT-WOPT method reaches it in significantly less time or with fewer iterations. Under~low-noise settings, the~(TT) algebraic initialization is already close to the optimal solution and thus does not require many iterations to converge. The~median number of iterations is indicated by vertical~bars.}
	\label{fig:expinit}
\end{figure}

A $4$th-order random tensor in TT format  \(\ten{x}\in \R^{40\times40\times40\times40}\)  is generated with a fixed  \(\mathrm{rank}_{TT}(\ten{x})=(1,5,5,5,1)\).  Random Gaussian i.i.d. noise is added to the generated tensor at SNR levels of 100 dB, 75 dB, 50 dB, and~25 dB to simulate low and moderate noise conditions. Sixty percent of the mode-4 fibers are randomly removed. Completion is then performed with two initialization strategies: one using 10 different random initializations, and~a hybrid approach where the TT approximation computed by the algebraic algorithm is used as the initialization for (one run of) TT-WOPT. Figure~\ref{fig:expinit} compares the number of iterations required to reach the same convergence criterion for  the algebraic and random initialization strategies. Each dot indicates the number of iterations required to reach the convergence criterion in a single experiment. For~the random strategy, each dot corresponds to the median number of iterations across 10 different initializations. A~total of 100 experiments is conducted. The~(median) accuracies achieved at SNR levels of 100 dB, 75 dB, 50 dB, and~25 dB are 
$4.51 \times 10^{-6}$, $3.27 \times 10^{-5}$, $2.54 \times 10^{-4}$, and~$2.76 \times 10^{-3}$, respectively, for~the hybrid strategy, and~
$1.28 \times 10^{-6}$, $3.01 \times 10^{-5}$, $2.52 \times 10^{-4}$, and~$2.76 \times 10^{-3}$, respectively, for~the random initialization strategy;  thus, both methods achieve comparable accuracies for the considered SNR levels. Meanwhile, with~algebraic initialization, the~TT-WOPT method reaches this accuracy in significantly fewer iterations, as~the algebraic initialization is already close to the optimal solution; see Figure~\ref{fig:expinit} (left).

\begin{figure}[H]
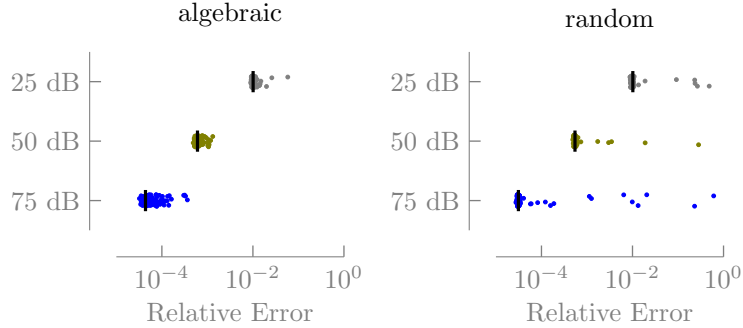

	\definecolor{mycolor2}{rgb}{0,0,1}
	\definecolor{mycolor3}{rgb}{0.502, 0.502, 0}
	\definecolor{mycolor4}{rgb}{0.5,0.5,0.5}
	
	\centering

	\caption{Performance comparison of TT-WOPT initialization methods at 25\,dB, 50\,dB, and~75\,dB SNR levels. While median relative errors (indicated by vertical bars) for random initialization ($0.0100$, $5.50 \times 10^{-4}$, and~$3.11 \times 10^{-5}$) are comparable to algebraic initialization ($0.0101$, $6.06 \times 10^{-4}$, and~$4.36 \times 10^{-5}$), respectively, random initialization exhibits significantly lower success rates in accurate TT completion. Algebraic initialization achieves reliable convergence in all trials, while random initialization fails to converge in a substantial fraction of cases, as~evidenced by the wider error~distributions.}
	\label{fig:expinitsA}
\end{figure}

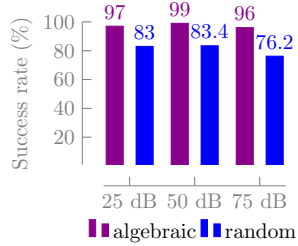
\begin{figure}[H]
	\definecolor{mycolor1}{rgb}{0.545,0,0.545}
	\definecolor{mycolor2}{rgb}{0,0,1}
	\centering
	\begin{tikzpicture}[scale=0.8]
		\begin{axis}[ width=0.35\textwidth,
			ybar,
			bar width=8pt,
			enlarge x limits=0.2,
			axis x line = bottom,
			axis y line =left,
			legend style={at={(0.5,-0.35)},
				anchor=north,legend columns=-1, draw=none},
			ylabel={Success rate (\%)},
			xtick={1,3,6},
			xticklabels={{25 dB},{50 dB},{75 dB}},
			xmin=0,
			xmax=7,
			ymin=1,
			ymax=100,
			xtick=data,
			nodes near coords,
			nodes near coords align={vertical},
			every axis x label/.append style={
				at={(axis description cs:0.5,-0.3)},
				anchor=north, color=gray
			},
			every axis y label/.append style={
				at={(axis description cs:-0.35,0.5)},
				anchor=south,
				rotate=0, color=gray},
			x axis line style={yshift=-10pt,-,color=gray, thin},
			xticklabel style={yshift=-10pt,color=gray},
			y axis line style={xshift=-10pt,-,color=gray},
			yticklabel style={xshift=-10pt,color=gray},
			xtick style={yshift=-10pt, line width=0.5pt},
			ytick style={xshift=-10pt, line width=0.5pt},
			xtick align=inside,
			ytick align=inside,
			]
			\addplot+[ybar, draw = mycolor1, fill=mycolor1, every node near coord/.append style={color=mycolor1}, bar shift=-7pt,] coordinates {(1,97) (3.5,99) (6,96)};
			\addplot+[ybar, draw = mycolor2, fill=mycolor2, every node near coord/.append style={color=mycolor2}, bar shift=7pt,] coordinates {(1, 83) (3.5,83.40) (6, 76.2)};
			\legend{algebraic,random}
		\end{axis}
	\end{tikzpicture}
	\caption{Each noise level has its threshold value: a completion is deemed successful if the relative error is below 0.015, 0.00125, and~0.000175 for SNR values of 25 dB, 50 dB, and~75 dB, respectively. The~algebraic strategy exhibits a higher success rate than the average success rate of the random strategy. This is also evident in Figure~\ref{fig:expinitsA}, where the standard deviation of the relative errors is lower for the algebraic strategy, indicating more consistent~performance.}
	\label{fig:expinitsB}
\end{figure}

Next,  we investigate the success of the two approaches in accurately estimating the underlying decomposition when the rate of missing fibers is high. A~$4$th-order random tensor in TT format  \(\ten{x}\in \R^{20 \times 20 \times 20 \times 20}\)  is generated with a fixed  \(\mathrm{rank}_{TT}(\ten{x})=(1,4,4,4,1)\).  Random Gaussian i.i.d. noise is added to the generated tensor at SNR levels of 75\,dB, 50\,dB, and~25\,dB. Approximately 20\% of the mode-4 fibers are sampled in such a way that the working conditions, as~discussed in Section~\ref{subsec:uniq}, on~the observation pattern are satisfied. Completion of the tensor is achieved similarly to the method described in the previous experiment. Figure~\ref{fig:expinitsA} shows the distribution of relative error across three different noise levels. Each dot represents one of 100 experiments, with~medians calculated over 10 initializations for the random strategy. It is observed that the TT-WOPT becomes highly successful when initialized with the algebraic method. To~assess the accuracy, we define thresholds for SNR values of 25\,dB, 50\,dB, and~75\,dB as 0.015, 0.00125, and~0.000175, respectively. A~completion is considered  successful if the relative error is below the threshold associated with the specified noise level. Figure~\ref{fig:expinitsB} shows the success rate at each noise level, defined as the proportion of 100 experiments that satisfy the accuracy criterion. The~success rate of the algebraic strategy is defined as the fraction of successful trials, whereas the success rate of the random strategy is computed by first calculating the fraction of successful trials for each of the 10 initializations and then averaging these fractions to obtain the final success rate. Note that each noise level has its corresponding threshold, meaning comparisons are made within the same SNR level rather than across different levels.  At~an SNR of 75 dB, the~algebraic strategy successfully completed 96 out of 100 experiments. In~comparison, the~random strategy has an average success rate of 76.2 out of 100 over 10 initializations. Overall, the~algebraic strategy consistently demonstrated a higher success rate than the random strategy under different noise~conditions.

\unskip

\unskip

\subsubsection{Proxy for Non-Negative~CPD}
\label{proxyCPD}

In this  experiment, we focus on computing a non-negative CPD of a non-negative data tensor that is observed through mode-$N$ fibers. The~constrained CPD is computed using two strategies: the \emph{full} strategy, where partial data in a dense format is directly used to compute the non-negative CPD, and~the \emph{compressed} strategy, where we first compute the TT approximation using our method and then use the compressed representation to compute the non-negative CPD. We use  projected  Gauss--Newton algorithm to compute the CPD (\texttt{cpd\_nls} with \texttt{nlsb\_gndl} solver), as~discussed in~\cite{nico2019compactrep}. 

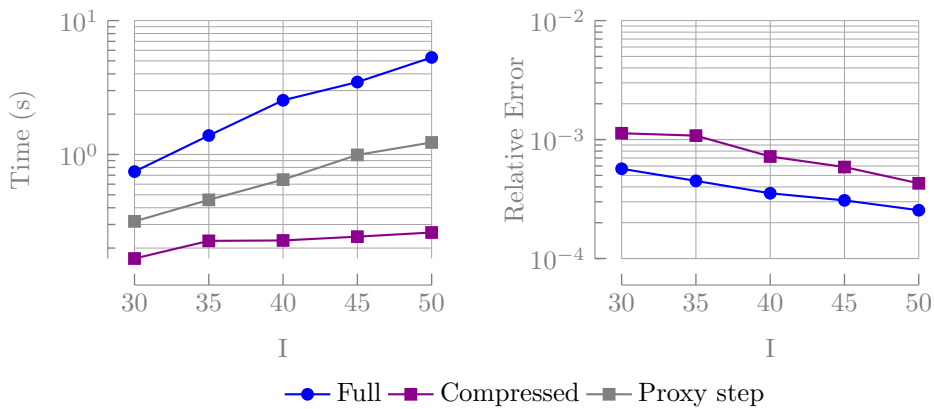
\begin{figure}[H]
	\centering
	\begin{tikzpicture}
		\begin{groupplot}[
			group style={group size=2 by 1, horizontal sep=2.5cm},
			width=0.42\textwidth,
			]
			\nextgroupplot[
			legend to name=sharedlegendexpcpd,
			legend style={draw=none, fill=none},
			legend columns=3,
			xlabel={I},
			ylabel={Time (s)},
			axis lines=middle,
			enlargelimits=true,
			axis x line = bottom,
			axis y line =left,
			grid=both,
			grid style={line width=0.2pt, draw=gray!70},
			ymax=10,
			xtick={30, 35, 40, 45,  50},
			ymode=log,
			every axis x label/.append style={
				at={(axis description cs:0.5,-0.3)},
				anchor=north, color=gray
			},
			every axis y label/.append style={
				at={(axis description cs:-0.3,0.5)},
				anchor=south,
				rotate=90, color=gray},
			x axis line style={yshift=-10pt,-,color=gray, thin},
			xticklabel style={yshift=-10pt,color=gray},
			y axis line style={xshift=-10pt,-,color=gray},
			yticklabel style={xshift=-10pt,color=gray},
			xtick style={yshift=-10pt, line width=0.5pt},
			ytick style={xshift=-10pt, line width=0.5pt},
			xtick align=inside,
			ytick align=inside,
			]
			\addplot [thick, blue, mark=*, mark size=2pt]
			table[row sep=crcr]{%
				30	0.74518755\\
				35	1.38566945\\
				40	2.5396397\\
				45	3.47358005\\
				50	5.30593935\\
			};
			\addlegendentry{Full}
			
			\addplot [thick, color={rgb,255:red,139; green,0; blue,139}, mark=square*, mark size=2pt]
			table[row sep=crcr]{%
				30	 0.16706665\\
				35	0.22632995\\
				40	0.2279663\\
				45	0.2433297\\
				50	0.2613543\\
			};
			\addlegendentry{Compressed}
			
			\addplot [thick, gray, mark=square*, mark size=2pt]
			table[row sep=crcr]{%
				30	0.3164918\\
				35	0.458194\\
				40	0.6488762\\
				45	0.9928255\\
				50	1.22989615\\
			};
			\addlegendentry{Proxy step}
			
			\nextgroupplot[
			xlabel={I},
			ylabel = {Relative Error},
			xtick={30, 35, 40, 45,  50},
			axis lines=middle,
			enlargelimits=true,
			axis x line = bottom,
			axis y line = left,
			ymode=log,
			ymax=0.01,
			ymin=0.0001,
			grid=both,
			grid style={line width=0.2pt, draw=gray!70},
			every axis x label/.append style={
				at={(axis description cs:0.5,-0.3)},
				anchor=north, color=gray
			},
			every axis y label/.append style={
				at={(axis description cs:-0.3,0.5)},
				anchor=south,
				rotate=90, color=gray},
			x axis line style={yshift=-10pt,-,color=gray, thin},
			xticklabel style={yshift=-10pt,gray=black},
			y axis line style={xshift=-10pt,-,color=gray},
			yticklabel style={xshift=-10pt,color=gray},
			xtick style={yshift=-10pt, line width=0.5pt},
			ytick style={xshift=-10pt, line width=0.5pt},
			xtick align=inside,
			ytick align=inside,
			]
			\addplot [thick, blue, mark=*, mark size=2pt]
			table[row sep=crcr]{%
				30	0.000568239790330996\\
				35	0.000449404332266409\\
				40	0.00035351466391419\\
				45	0.000308227569563586\\
				50	0.000254726579098751\\
			};
			
			\addplot [thick, color={rgb,255:red,139; green,0; blue,139}, mark=square*, mark size=2pt]
			table[row sep=crcr]{%
				30	0.00112973844565031\\
				35	0.00107941533680615\\
				40	0.000720449931423547\\
				45	0.000587609734645112\\
				50	0.000428921122874142\\
			};
			
		\end{groupplot}

		\node at ($ (group c1r1.south)!0.5!(group c2r1.south) + (0,-1.8cm) $){\pgfplotslegendfromname{sharedlegendexpcpd}};
		
	\end{tikzpicture}
	\caption{(\textbf{Left}) The time required to compute the non-negative CPD, given a TT approximation as a prior, increases very slowly with problem size. On~the other hand, the~time required to compute the (TT) proxy  increases, as~indicated by the grey curve. However, the~overall computation time in the compressed strategy remains significantly lower than that of the full strategy. (\textbf{Right}) The trade-off is that the compressed strategy is slightly less accurate, but~this can be resolved by refining the proxy for a few~iterations.}
	\label{fig:cpdexp}
\end{figure}

A 4th-order tensor in  CPD format $\ten{X} \in \R^{I \times I \times I \times I}$ is generated with a fixed  $\text{rank}_{CP}(\ten{x}) = 4 $ by sampling the factor matrices from a uniform distribution  $\ten{U}(0, 1)$.  Random Gaussian i.i.d.  noise is added to the generated tensor to achieve an SNR of 40 dB. Fifty percent of the mode-4 fibers are randomly removed from the noisy tensor. The~constrained rank-4 CPD is computed using both the full and compressed strategies for \( I \) varying from 30 to 50. The~median of the computation time is plotted across 50 trials in Figure~\ref{fig:cpdexp} (left). It is observed that the use of the TT approximation as a prior significantly improves speed. The~overall computation time in the compressed strategy, which includes the proxy step and the computation of the constrained CPD from the compressed tensor, remains lower than that of the full strategy. The~small accuracy gap that arises in noisy settings can be eliminated by refining the proxy over a few iterations. Moreover, the~computational cost of the proxy step can be further reduced by using parallel implementations or randomized SVD~algorithms.

\unskip

\subsection{Note on Accuracy Gap Between the Algebraic and Optimization Methods}

The observed accuracy gap between the algebraic and optimization methods stems not only from the algebraic method's reliance solely on standard NLA operations---which makes it slightly less accurate in noisy cases---but also from the difference in the proportion of data points that each method uses. When 70\% of the fibers along a specific mode are randomly removed, the~optimization method uses all of the remaining 30\% of the observations to compute the completion. In~contrast, the~algebraic method only uses the observed submatrices that satisfy the working conditions for recovering the column spaces of the matrix unfoldings, discarding the rest. Consequently, the~algebraic method uses fewer than 30\% of the observed fibers. Under~deterministic sampling, where all 30\% of the observed fibers satisfy the working conditions, the~accuracy gap is expected to be~small. \par

Another key point is that only the observed rows of the (individual) mode-2 slices, as~in Equation~(\ref{eqn:ttfw1}), are used in the above experiments. The~observed submatrices spanning multiple slices, satisfying the working conditions, can also be utilized. Let us consider an example of an observed submatrix spanning two mode-2 slices, $\mat{s}^{(l_1)}_{r}\tilde{\mat{x}}_{:l_1:} \in \mathbb{R}^{| \alpha_{l_1} | \times I_N}$ and $\mat{s}^{(l_2)}_{r}\tilde{\mat{x}}_{:l_2:} \in \mathbb{R}^{| \alpha_{l_2} | \times I_N}$, which have overlapping rows indexed by $\alpha_{l_1 \cap l_2} = \alpha_{l_1} \cap \alpha_{l_2}$, where $| \alpha_{l_1 \cap l_2}|  = K \geq R$. Then, the~concatenated matrix \(\mat{s}^{(l_1 \cap l_2)}_{r} \left[ \tilde{\mat{x}}_{:l_1:} \;\; \tilde{\mat{x}}_{:l_2:} \right] \in \mathbb{R}^{K \times 2I_N}
\), where the row selection matrix $\mat{s}^{(l_1 \cap l_2)}_{r}$ selects rows indexed by $\alpha_{l_1 \cap l_2}$, can also serve as an additional observed submatrix in the subspace computation approaches discussed in \mbox{Sections~\ref{subsec:ssc} and~\ref{subsec:ssI}}. As~a result, the~column subspaces can be estimated more accurately by incorporating these additional submatrices. However, this approach comes with additional computational cost as it requires identifying the overlapping slices. Computing the SVDs of these larger matrices can also be relatively expensive, but~the results are expected to be close to those obtained by optimization methods. In~the remainder of this section, we demonstrate this through an experiment. We compute the (TT) approximation of an incomplete tensor observed fiber-wise using the algebraic algorithm with two different approaches: one that utilizes only individual slices and another that incorporates both individual mode-2 slices and the largest observed submatrix from each pair of~slices. \par
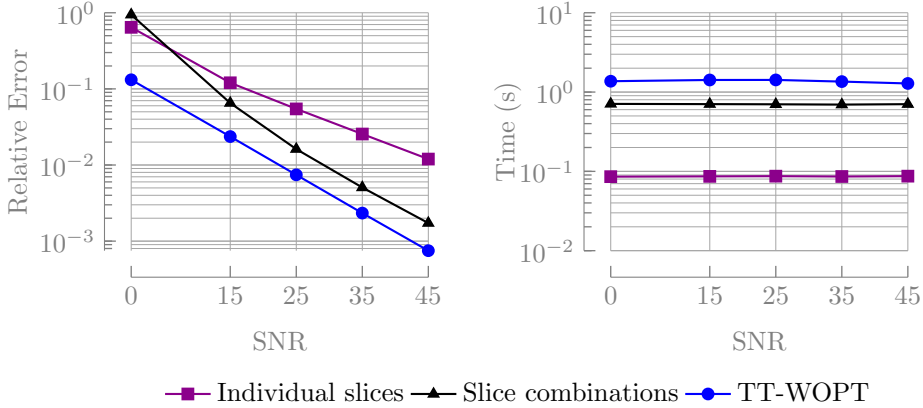
\begin{figure}[H]
	\begin{tikzpicture}
		\begin{groupplot}[
			group style={group size=2 by 1, horizontal sep=2.4cm},
			width=0.42\textwidth,
			legend to name=sharedlegendepxlast,
			legend style={draw=none},
			legend columns=3,
			]
			\nextgroupplot[
			xlabel={SNR},
			ylabel={Relative Error},
			axis lines=middle,
			enlargelimits=true,
			axis x line = bottom,
			axis y line =left,
			grid=both,
			grid style={line width=0.2pt, draw=gray!70},
			xtick={0,15,25,35,45},
			ymode=log,
			ymin=0,
			ymax=1,
			scale=1,
			every axis x label/.append style={
				at={(axis description cs:0.5,-0.3)},
				anchor=north, color=gray
			},
			every axis y label/.append style={
				at={(axis description cs:-0.32,0.5)},
				anchor=south,
				rotate=90, color=gray},
			x axis line style={yshift=-10pt,-,color=gray},
			xticklabel style={yshift=-10pt,color=gray},
			y axis line style={xshift=-10pt,-,color=gray},
			yticklabel style={xshift=-10pt,color=gray},
			xtick style={yshift=-10pt, line width=0.5pt},
			ytick style={xshift=-10pt, line width=0.5pt},
			xtick align=inside,
			ytick align=inside,
			]
			
			\addplot[thick, color={rgb,255:red,139; green,0; blue,139}, mark=square*, mark size=2pt]  
			table[row sep=crcr]{%
				0	0.643540055833742\\
				15	0.120229228916706\\
				25	0.0545671681642027\\
				35	0.0255270012650784\\
				45	0.0119686743483476\\
			};
			\addlegendentry{Individual slices}
			
			\addplot[thick, black, mark=triangle*, mark size=2pt] 
			table[row sep=crcr]{%
				0	0.942443378611012\\
				15	0.0651372154445616\\
				25	0.0161728311629729\\
				35	0.00506689734808574\\
				45	0.00172647633855791\\
			};
			\addlegendentry{Slice combinations}
			
			\addplot[thick, color=blue, mark=*, mark size=2pt]  
			table[row sep=crcr]{%
				0	0.131675494537895\\
				15	0.0236391638697714\\
				25	0.00742370685018181\\
				35	0.00233558950505461\\
				45	0.000749106429390686\\
			};
			\addlegendentry{TT-WOPT}
			
			\nextgroupplot[
			xlabel={SNR},
			xtick={0,15,25,35,45},
			ylabel={\textcolor{gray}{Time (s)}},
			axis lines=middle,
			enlargelimits=true,
			axis x line = bottom,
			axis y line = left,
			ymode=log,
			ymin=0.01,
			ymax=10,
			grid=both,
			grid style={line width=0.2pt, draw=gray!70},
			every axis x label/.append style={
				at={(axis description cs:0.5,-0.3)},
				anchor=north, color=gray
			},
			every axis y label/.append style={
				at={(axis description cs:-0.27,0.5)},
				anchor=south,
				rotate=90, color=gray},
			x axis line style={yshift=-10pt,-,color=gray},
			xticklabel style={yshift=-10pt,color=gray},
			y axis line style={xshift=-10pt,-,color=gray},
			yticklabel style={xshift=-10pt,color=gray},
			xtick style={yshift=-10pt, line width=0.5pt},
			ytick style={xshift=-10pt, line width=0.5pt},
			xtick align=inside,
			ytick align=inside,
			]

			\addplot[thick, color={rgb,255:red,139; green,0; blue,139}, mark=square*, mark size=2pt]  
			table[row sep=crcr]{%
				0	0.0857826\\
				15	0.0864826\\
				25	0.087001\\
				35	0.0861555\\
				45	0.0870402\\
			};
			\addlegendentry{Individual slices}
			
			\addplot[thick, black, mark=triangle*, mark size=2pt]
			table[row sep=crcr]{%
				0	0.70933775\\
				15	0.7046453\\
				25	0.7009625\\
				35	0.6947237\\
				45	0.70378575\\
			};
			\addlegendentry{Slice combinations}
			
			\addplot[thick, color=blue, mark=*, mark size=2pt] 
			table[row sep=crcr]{%
				0	1.37010755\\
				15	1.42024945\\
				25	1.42257025\\
				35	1.35547355\\
				45	1.28402365\\
			};
			\addlegendentry{TT-WOPT}
			
		\end{groupplot}

		\node at ($ (group c1r1.south)!0.5!(group c2r1.south) + (0,-1.9cm) $) 
		{\pgfplotslegendfromname{sharedlegendepxlast}};
	\end{tikzpicture}
	\caption{(\textbf{Left}) Accuracy of the algebraic method increases as observed submatrices spanning pairs of slices are included in the computation of the column spaces of the matrix unfoldings, and~it becomes close to that of TT-WOPT. (\textbf{Right}) However, this improvement comes with additional computational~cost.}
	\label{fig:expslicecomb}
\end{figure}

A 4th-order tensor in TT format  \(\ten{x} \in \R^{16 \times 16 \times 16 \times 16}\) is generated with a fixed  \(\text{rank}_{TT}(\ten{x}) = (1, 3, 3, 4, 1)\). Random Gaussian i.i.d. noise is then added to the generated tensor with SNR varying from 0 dB to 45 dB. Sixty percent of the mode-4 fibers are randomly removed from the noisy tensor. Completion is then performed using the algebraic algorithm with the two strategies. Figure~\ref{fig:expslicecomb} shows that using the additional observed submatrices significantly improves accuracy, bringing it close to that of TT-WOPT. However, this approach is slightly slower due to the extra computational overhead, as mentioned previously. Overall, the~computational cost remains low compared to TT-WOPT, while the accuracy gap is very small. It is worth noting that, in~this experiment, we only combined pairs of slices; however, it may be possible to achieve an accuracy closer to TT-WOPT by selecting triplets or more~slices.

\section{Conclusions}\label{sec:con}
We introduced an algebraic framework for computing the TT decomposition of an incomplete tensor observed fiber-wise (along a single specific mode). This framework builds on established algebraic techniques known for CPD and MLSVD~\cite{mikael2019fibersamp, stijn2023mlsvdfsj}. The~TT decomposition combines the features of MLSVD and CPD, making the extension of such algebraic methods in TT format highly consequential. The~proposed approach relies solely on standard NLA operations and is fast, being guaranteed to work under reasonable deterministic conditions on the observation~pattern.

We provided theoretical insights into piecewise subspace learning, an~essential ingredient of our method, and~discussed both the algebraic and generic uniqueness conditions for retrieving the TT cores (up to basis transformation).

Convincing numerical experiments demonstrate that our proposed approach is practical and valuable for real-life applications. When compared with state-of-the-art methods, our method has been observed to be fast, in~line with expectations, while achieving competitive accuracy in recovering partially (fiber-wise) observed tensors. Moreover, our experiments show that the solution obtained from the proposed method can also serve as a proxy for efficient subsequent computations, including the initialization of optimization-based methods, which are typically more computationally~expensive. 

\vspace{6pt}

   \section*{\large Declarations}
   \begin{itemize}
   	\item \textbf{Author Contributions:} Conceptualization, L.D.L. and S.S.S.; methodology, L.D.L. and S.S.S.; software, S.S.S.; validation, S.S.S.; formal analysis, S.S.S.; resources, L.D.L.; writing—original draft preparation, L.D.L. and S.S.S.;  project administration, L.D.L.; funding acquisition, L.D.L. All authors have read and agreed to the published version of the~manuscript.
   	
   	\item \textbf{Funding:} This work was supported by the Flemish Government's AI Research Program and KU Leuven Internal Funds (iBOF/23/064, C14/22/096, IDN/19/014). Lieven De Lathauwer and Shakir Showkat Sofi are affiliated with Leuven.AI - KU Leuven institute for AI, B-3000, Leuven, Belgium.

   	\item \textbf{Code availability:}  The source code for the proposed algorithm, along with all data necessary to reproduce the results, will be made available at \url{https://www.tensorlabplus.net/}. 
   	
   	\item \textbf{Conflict of interest:}  The funders had no role in the design of the study; in the collection, analyses, or~interpretation of data; in the writing of the manuscript; or in the decision to publish the~results.
   \end{itemize}
   
   \section*{\large Abbreviations}
   	The following abbreviations are used in this manuscript:	\\
   	
   	\noindent 
   	\begin{tabular}{@{}ll}
   		CPD & Canonical Polyadic Decomposition\\
   		MHR & Multidimensional Harmonic Retrieval\\
   		MLSVD & Multilinear Singular Value Decomposition\\
   		NLA & Numerical Linear Algebra\\
   		NP & Nondeterministic Polynomial Time\\
   		RMSE &  Root Mean Square Error\\
   		SiLRTC-TT & Simple Low-Rank Tensor Completion via Tensor Train\\
   		SNR & Signal-to-Noise Ratio\\
   		TMAX & Maximum Temperature\\
   		TMac-TT & Parallel Matrix Factorization for Low-Rank Completion via Tensor Train\\
   		TT & Tensor Train\\
   		TT-SVD & Tensor Train Singular Value Decomposition\\
   		TT-WOPT & Tensor Train Weighted Optimization\\
   	\end{tabular}

\end{document}